\documentclass[aip,graphicx]{revtex4-1}
\usepackage{amsmath}
\usepackage{amssymb}
\usepackage{natbib}
\newtheorem{notation}{Notation}[section]
\newtheorem{theorem}{Theorem}[section]
\newtheorem{proposition}{Proposition}[section]
\newtheorem{remark}{Remark}[section]

\newtheorem{definition}{Definition}[section]
\newtheorem{corollary}{Corollary}[section]
\draft 
\begin{document}
\title{Conformal Scale Geometry of Spacetime -- A lower bound  for a total mass} 



\author{Jean Baptiste PATENOU}
\email[]{jeanbaptiste.patenou@univ-dschang.org,\;patenou@gmail.com}
\affiliation{University of Dschang, Cameroon; Faculty
of Science, Department of Mathematics and Computer Science.
 Research Unit of Mathematics and Applications (RUMA),
              }

\date{\today}
\begin{abstract}
We devise a new approach for the analysis of issues of geometric pathologies
and black holes of a spacetime, based on a new mass function defined on an
ideal un-physical spacetime which models time-flow or time dilation.
The mass function is interpreted as an "extra" local energy density
that encodes the rate at which time comes to a "stop" (hardly
visible) or it measures how quickly the (illusory) Event horizon forms. This
latter is defined on the manifold with corners resulting from an
appropriate conformal compactification of the original physical
space-time, the concrete choice of compactification being tied to
the geometric structure of collapsing spacetimes. We define the
(illusory) Event horizon as the set of zero-mass function and
provide conditions for which it stands as a "black hole's event
horizon". As a first main result owing to the new definitions here,
we establish the existence of a lower bound for the "total mass", provided some
conditions on the extrinsic curvature of the space-time are
satisfied. The proof builds on geometric flow techniques. Namely, by
flowing the "black hole's event horizon", one is able to derive via a
Lagrangian formulation, a Minimization Problem for the "total mass"
and which is addressed under Iso-perimetric constraints' perspective. Prior to this main result, we first
provide hypotheses which assure the trapped surface's formation in
this context.

\end{abstract}


\maketitle 

\section{Introduction -- Main results}
\label{intro} One fundamental avenue in understanding the Einstein's
theory of general relativity and its implications is the study of
issues of singularities and black holes, which are fascinating
objects predicted by this theory
\cite{11}, \cite{17}, \cite{42},\cite{51},\cite{67},\cite{77}. Indeed,
several efforts have been made in the analysis of the formation of
singularities and black holes, the study of their nature and their
stability \cite{11}-\cite{18},\cite{42},
\cite{27},\cite{36},\cite{49},\cite{69}-\cite{77},\cite{82}. Some
remarkable results for stability of blackholes include: - the proof
of global nonlinear stability of Minkowski spacetime
\cite{4},\cite{18},\cite{45},\cite{53},\cite{54},\cite{80}, - the
proof of global nonlinear stability of Schwarzschild spacetime under
polarized perturbations \cite{94}, - the proof of global nonlinear
stability of the Kerr-de Sitter family of blackholes \cite{93}.
However, many facets of the theory of general relativity remain
questionable. As instance, the question of stability of the Kerr
black holes, or the full proof of stability conjecture still open
\cite{95}. Importantly, some open questions related to issues of
singularities and black holes are: - the definition of black hole as
a feature of spacetime, - the question of existence or definition of
horizons \cite{13},\cite{25},\cite{56},
\cite{65},\cite{66},\cite{75},\cite{78},\cite{91},  - their nature
and their stability--the fate of the Cauchy horizons - the
deterministic character of the theory of general relativity
\cite{2},\cite{3},\cite{10}-\cite{16},\cite{22}-\cite{37},\cite{82}.
These problems harbor two conjectures known as weak and strong
cosmic censorship conjectures formulated initially by Roger Penrose
\cite{66}. These conjectures have received many formulations, have
been to date at the core of many other scientific works in the field
of general relativity and still power many research perspectives
including, the question of definition of a general appropriate
quasi-local mass in general relativity
\cite{20},\cite{21},\cite{85}-\cite{90},\cite{92}, the positivity of
the total mass or the question of establishment of a positive lower
bound for the total energy or mass of isolated systems in general
resumed as Penrose's conjecture. The Penroses's Conjecture asserts
that the total mass of a spacetime including Black holes of area $A$
is at least $\sqrt{\frac{A}{16\pi}}$. It dates to $1973$ \cite{65}.
This beautiful lower bound inequality has been the source of many
interesting research proposals in recent years leading to some
important results such as: the proof of the positive mass theorem,
the proof of the Riemannian Penrose-Inequality, and more recently
the proof of the Null Penrose Conjecture
\cite{1},\cite{5}-\cite{9},\cite{40},\cite{47},\cite{48},
\cite{55},\cite{57},\cite{60},\cite{70}-\cite{72},\cite{74},\cite{79},\cite{90}.
There are other questions raised by quantum gravity which make the
satisfaction not yet within the reach for scientists of general
relativity, these include in particular the definition of a horizon
and its existence, and the hypothetical phenomenon of "Firewall". These reasons justify somewhat
the emergence of theories of modified gravity.\\
Efforts in the comprehension of the general theory of
relativity and its implication, and progress obtained for the issues
mentioned above rest basically on the followings:\\ -- mastery of
advanced tools of differential geometry and in particular the study
of Riemannian and Lorentzian structures
\cite{11},\cite{17},\cite{19},\cite{42},\cite{51},\cite{77},
\cite{81},\\-- study of exact solutions of the Einstein equations
(Minkowski, Schwarzschild, Kerr, Reissner-Norstr$\ddot{o}$m,... )
notably their stability
\cite{4},\cite{18},\cite{45},\cite{53},\cite{54},\cite{80},
\cite{27},\cite{36},\cite{49},\cite{76},\cite{77},\cite{82},\cite{93},\cite{94},\\
-- systematic construction of spacetimes by solving Cauchy problems
and the analysis of properties of the obtained solutions (initial
data constraints problem, evolution problem, existence theorems for
differential equations, asymptotic properties,...)
\cite{11},\cite{15},\cite{38},\cite{42},\cite{50},\cite{51},\cite{62},
\cite{63},\\ -- analysis of waves on fixed backgrounds (decay
properties, boundedness,...)
\cite{2},\cite{3},\cite{10}-\cite{16},\cite{22}-\cite{37},\cite{39},\cite{41}-\cite{44},
\cite{46},\cite{48}-\cite{52},\cite{58}-\cite{61},\cite{82}-\cite{84},\\
- and even numerical relativity.\\ Concerning methods of studying
asymptotic questions in general relativity, many authors agree in
the conformal treatment of infinity as conceived by Penrose
\cite{64},\cite{68}, and this substantiates also our present
analysis.\\
 This research paper which is a revised version of the first
  part of my preprint
 \cite{98}, aims at contributing to some of the
issues above. For this purpose, we adopt a mass function  modeling
the time-flow or the time dilation, or which helps to encode
time-visibility at every point, defined on an appropriate conformal
manifold (with corners) associated to the original physical
spacetime. This unphysical conformal manifold is obtained by
conformally embedding the spacetime under consideration in a new
manifold with boundaries via a specific gauge or coordinates system.
Given that the current mass function determines positions in the
road of time, we localize the "blackhole's event horizon" candidate
$\mathcal{H}$ as the hypersurface of zero-mass meaning that time is
"stopped" (hardly visible), that is the intrinsic geometry of the
Event horizon is required to be time independent, whereas the
geometry outside may be dynamical and may admit gravitational or
other radiations. Using the wave fronts sets, the Event horizon can
be scended in an open subset $\mathcal{H}^+$ designed as the "actual
event horizon" and a subset $\mathcal{H}^-$ which we refer to as the
"Apparent horizon". This new approach postulates a possible
cohabitation of the (possible stationary) causal future of
$\mathcal{H}^+$ (i.e. "the Black hole's region") and the causal future
of $\mathcal{H}^-$ expected to be the White hole and likely to be
instable. It appears also in this setting that any mass
(concentration or blow-up) singularity is "preceded" by an illusory
event horizon (which may coincide with the black hole's event
horizon) giving an approach to the formulation of the Weak Cosmic
Censorship Conjecture. We signal that the exact relation between the
mass inflation here and singularities as incompleteness of causal's
geodesics remains a point to examine.
 We also emphasize as one can remark that the standard definition of black hole based on the
idea of a global event horizon is abandoned here in favor of a new
one which relies on "time hardly visible in an appropriate gauge at
later times".\\
The main results of this paper are twofold and consist in a proof of
a theorem relative to the formation of trapped surfaces, and a
proof of the existence of a lower bound  for  the "total mass" provided some hypotheses are given for the isoperimetric constraint. These results are achieved according to
the new mass function, the related quantities and sets, and by
formulating (thanks to geometric flow techniques and via a Lagrangian formulation) the lower bound
problem for the "total mass" as a Minimization problem under a dynamical
Isoperimetric constraint \cite{81}.  We emphasize that in the solution process the dynamical integrand in the "total mass" satisfies the Euler-Lagrange equation, this offers rich perspective in terms of the choice of the Ispoerimetric constraint -- it might involve specific: energies, area or volume, ...; also the flow techniques might involve the flow of the metrics and/or the flow of a particular surface depending on the particular or desirable needs or ends.
{\theorem {Let $(V,g)$ be a global
in time
 spacetime solution of the Einstein field
 equations where $g$ is of the following form in coordinates $(x^\alpha)$:
\begin{equation}
    g=-\frac{|\overline{g}|}{\Omega^2}(dx^0)^2+\overline{g}_{ij}dx^idx^j,\;\overline{g}=(\overline{g}_{ij}),
\end{equation} $|\overline{g}|$ is the determinant of $\overline{g}$, $\Omega\equiv \Omega (x^i)$ is an arbitrary positive
given scalar density satisfying $\Omega\rightarrow 0$ as
$r=|(x^i)|\rightarrow +\infty$. \\ Let denote $h=\Omega^2 g$
the conformal metric, $G$ its inverse, and $\widehat{\nabla}$ its associated riemannian
connection. Let be specified the gauge $(\omega^\alpha)$ for the spacetime $(V,g)$, with $
\omega^0=e^{-\Omega x^0},\;\omega^1=\Omega,\;\omega^a=\arctan
(\Omega x^a)$.\\ One suppose that $\mathcal{H}=\{m\equiv G\left(\frac{d\omega^0}{\Omega\omega^0},\frac{d\omega^0}{\Omega\omega^0}\right)=0\}\neq \emptyset$ and defines a (smooth) hypersurface, and the attached
 geometric objects (in $(\widehat{V},h)$), $\textbf{U},\textbf{V}$ and quantities $\textbf{a},\textbf{b}^\pm$ as described in
       (\ref{31})-(\ref{33}) satisfy:
  \begin{equation}
    h\left(\widehat{\nabla}_{\frac{\partial}{\partial\omega^0}}^{\frac{\widehat{\nabla}\omega^0}{\Omega\omega^0}},
\frac{\widehat{\nabla}\omega^0}{\Omega\omega^0}\right) <0, \;
   h\left(\widehat{\nabla}_{\widehat{h}^{0i}\frac{\partial}{\partial\omega^i}}^{\frac{\widehat{\nabla}\omega^0}{\Omega\omega^0}},
\frac{\widehat{\nabla}\omega^0}{\Omega\omega^0}\right)>0;
 \end{equation}
 \begin{equation}
  tr\textbf{V}<0,\;
      tr\textbf{U}\in\left]-\frac{tr\textbf{V}}{\textbf{b}^-}
      ,-\frac{tr\textbf{V}}{\textbf{b}^+}\right[.
       \end{equation}Then the spacetime $(V,g)$ undergoes trapped surfaces's formation.  }\label{th5}}
{\theorem {One supposes that a gauge $(\omega^\alpha)$ is specified as in theorem \ref{th5} for a spacetime $(V,g)$, $\mathcal{H}=\{m\equiv G\left(\frac{d\omega^0}{\Omega\omega^0},\frac{d\omega^0}{\Omega\omega^0}\right)=0\}\neq \emptyset$ and defines a (smooth) hypersurface,
  and the attached
 geometric objects (in $(\widehat{V},h)$) $\textbf{U},\textbf{V}$, and quantities $\textbf{a},\textbf{b}^\pm$ as described in
       (\ref{31})-(\ref{33}) satisfy:
  \begin{equation}
    h\left(\widehat{\nabla}_{\frac{\partial}{\partial\omega^0}}^{\frac{\widehat{\nabla}\omega^0}{\Omega\omega^0}},
\frac{\widehat{\nabla}\omega^0}{\Omega\omega^0}\right) <0, \;
   h\left(\widehat{\nabla}_{\widehat{h}^{0i}\frac{\partial}{\partial\omega^i}}^{\frac{\widehat{\nabla}\omega^0}{\Omega\omega^0}},
\frac{\widehat{\nabla}\omega^0}{\Omega\omega^0}\right)>0;
 \end{equation}
 \begin{equation}
  tr\textbf{V}<0,\;
      tr\textbf{U}\in\left]-\frac{tr\textbf{V}}{\textbf{b}^-}
      ,-\frac{tr\textbf{V}}{\textbf{b}^+}\right[.
       \end{equation} \\One supposes further that for a given smooth function $P_0$  of three real variables, the fundamental equation (\ref{e1}) (with $P$ described in (\ref{e3}))
              \begin{equation*}
    \frac{\partial P_0}{\partial
\widetilde{T}}-\frac{d}{ds}\left(\frac{\partial P_0}{\partial
    \widetilde{T}'}\right)+\left[\frac{1}{\lambda}\frac{\partial P}{\partial \widetilde{T}'}(\widetilde{T}(+\epsilon),\widetilde{T}'(+\epsilon))
+ \frac{\partial P_0}{\partial \widetilde{T}'}(\widetilde{T}(+\epsilon),\widetilde{T}'(+\epsilon))\right]\delta_{+\epsilon}-
    \end{equation*}
\begin{equation}\label{e2}
\left[\frac{1}{\lambda}\frac{\partial P}{\partial \widetilde{T}'}(\widetilde{T}(-\epsilon),\widetilde{T}'(-\epsilon))
+ \frac{\partial P_0}{\partial \widetilde{T}'}(\widetilde{T}(-\epsilon),\widetilde{T}'(-\epsilon))\right]\delta_{-\epsilon}
   =0
\end{equation}
admits a solution $\widetilde{T}_\ast (.,\lambda)$ with $\lambda$ satisfying $I(\widetilde{T}_\ast)=K$.\\ Then the total mass $\mathfrak{M}$ corresponding to
the mass function
$m=G(\frac{d\omega^0}{\omega^1\omega^0},\frac{d\omega^0}{\omega^1\omega^0})\equiv
h(\frac{\widehat{ \nabla}\omega^0}{\omega^1\omega^0},\frac{\widehat{\nabla}\omega^0}{\omega^1\omega^0})$
(with $G$ the inverse of $h\equiv \Omega^2g $) admits a lower bound .  }\label{th6}}\\

 The different steps of the paper comprise:\\
 - the preliminaries related to the Cauchy problem for the Einstein
 equations and the role of gauges,-- the choice of temporal gauge offers the means of exhibiting a
 conformal factor for a partial (radial) compactification of spacetime,\\
 - the compactification of spacetime using appropriate conformal
 factor,\\
 - the definitions of geometric quantities and sets, and a derived geometric analysis,
  one establishes some hypotheses relative to trapped surfaces's formation,\\
 - the description and resolution of a lower bound problem for the total mass
 as a minimization problem with an associated Isoperimetric
 constraint,
 using a foliation of the compactified spacetime by
 a deformation of the "Black hole's event horizon", this yields one fundamental equation with important perspectives
  depending on the choice of the isoperimetric constraint,\\
- the conclusion and outlook.
\section{Preliminaries}
The Einstein's theory of general relativity postulates a spacetime
$(\mathcal{M},g)$ which must satisfy the Einstein's equations (in
units where $c=G=1$)
\begin{equation}\label{}
    Ric (g)-\frac{1}{2}R_{scal}\; g=8\pi T,
\end{equation}where $Ric (g)$ is the Ricci curvature tensor of $g$, $R_{scal}$
the corresponding scalar curvature, and $T$ the stress
energy-momentum tensor of matter.
\subsection{Cauchy problem--Gauge}
 A systematic way to construct
a spacetime without any symmetry assumption is by solving a Cauchy
problem for the Einstein equations where initial data are prescribed
on an appropriate initial hypersurface
\cite{11},\cite{42},\cite{51},\cite{77}.
 The Einstein's equations-- geometric in nature, however read as a set
of ill-posed complicated system of partial differential equations in
arbitrary gauge or coordinates system
\cite{11},\cite{42},\cite{51},\cite{77} for the unknown metric. It
is usual to require that the coordinate system satisfy
specific properties and this induces a choice of gauge. One known
gauge is Harmonic gauge or wave gauge with its generalization the
generalized wave map gauge used by P. T. Chru\`sciel et al.
(see \cite{21} and references in). Other known gauges are the Double
null foliation gauge \cite{51}, and the Temporal gauge \cite{11}. The
gauge has the role of splitting the Einstein's equations to an
evolution system and the set of equations known as the constraint's
equations. Two types of problems then arise which are the initial
data constraints's problem and the problem of evolution of initial
data. The first consists of studying how to effectively prescribe the
initial data satisfying the constraints's equations on the initial
manifold for the evolution system obtained according to the gauge,
so that its solution yields the solution of the full
 system of Einstein's equations. The initial manifold may be a spacelike
  hypersurface (spacelike Cauchy problem) or
one or more null hypersurfaces (characteristic Cauchy problem). In
the case of a spacelike Cauchy problem, the set of constraints are
standard \cite{11} whereas in the case of
characteristic hypersurfaces, they include standard constraints and
other gauge-dependent constraints (see
\cite{21},\cite{38},\cite{62},\cite{63} and references in). In the
characteristic Cauchy problem setting, the constraints are a set of
propagations equations along null geodesics generating the initial
hypersurfaces, so that the main difficulties in the resolution of
the initial data constraint's problem in this setting reside rather
in the null geometric description and construction of such
constraints (see \cite{21},\cite{38},\cite{62},\cite{63} and
references in).
\subsection {Trouble with gauges}
The gauge has the merit of providing an evolution system of partial
differential equations which can be solved by standard methods.
However, there is a trouble in choosing the right one since the
universe does not hand us with a preferred coordinates system which
always works. There is in fact a possibility of drawing wrong
conclusions about the properties of physical space because the
coordinates system is inadequate. This problem is at the heart of
issues of singularities and black holes (Measurement of Decay of
waves, Stability,...), -- new results of Klainermann--Szeftel
\cite{94},\cite{95}; and Hintz--Vasy \cite{93}, on the stability problem
involve new techniques for finding the right coordinates system.
Let's recall that the equivalence principle in General Relativity
(GR) is equivalent to an accelerating reference frame.
\subsection{Temporal gauge}
The gauge under consideration here is temporal gauge. If
$(x^\alpha)$ is a system of coordinates of the spacetime $(V,g)$,
the temporal gauge condition for this system reads
\begin{equation}
    \nabla^\lambda\nabla_\lambda x^0=0,\;g_{0i}=0;
\end{equation}$\nabla^\lambda\nabla_\lambda $ is the wave operator attached to $g$ --  $\nabla$ is the
covariant derivative operator attached to the metric $g$),; this induces
a metric of the form
\begin{equation}\label{1}
g=-\tau^2
(dx^0)^2+\overline{g}_{ij}dx^idx^j,\;\tau=c(x^i)\sqrt{|\overline{g}|},
\end{equation}where the function $(x^i)\rightarrow c(x^i)$ is an
arbitrary scalar density on the slices $x^0=t$,
$\overline{g}=(g_{ij})$ is the induced Riemannian metric on
$x^0=t,\;t\in \mathbb{R}$ \cite{11}. The evolution system attached
to the Einstein's equations in temporal gauge reads \cite{11}
\begin{equation}\label{2}
 \Omega_{ij}\equiv   \partial_0(R_{ij}-\Lambda_{ij})-\overline{\nabla}_i(R_{0j}-\Lambda_{0j})
    -\overline{\nabla}_j(R_{i0}-\Lambda_{i0})=
    0,
\end{equation}where
$\Lambda_{\mu\nu}=T_{\mu\nu}+\frac{T^\lambda_\lambda}{1-n}g_{\mu\nu}$.
For $c>0$, and for fixed $\overline{g}, \;T$, the system
$\Omega_{ij}=0$ is a quasi-diagonal system
 of wave equations for the extrinsic curvature $K=(K_{ij})$ of the
hypersurface $x^0=t$ \cite{11}. Together with the system
\begin{equation}\label{3}
    \partial_0\overline{g}_{ij}=-2\tau K_{ij},
\end{equation}they form a hyperbolic Leray system for $\overline{g}$
and $K$ provided $T$ is fixed, otherwise it is coupled with an
appropriate system of equations describing the matter contents of
the involved system. In reference \cite{38}, by solving the characteristic
initial data constraints problem related to such a system for
Vlasov-scalar field matter, it appeared that the scalar density $c$
is of specific form $c(x^i)=\frac{1}{|\gamma|}$; is obtained
actually by a null geometric description and resolution of the
initial data constraints, and does not depend on the matter contents
but only on a part of the free initial gravitational data. As our
attention is focused here on the causal features of the expected
global in time spacetime, this latter serves as a natural conformal
factor for a partial (radial) "bordification" of the spacetime; $c$
is therefore the boundary defining function for
infinity. Another reason in choosing the temporal gauge here
corresponds to our intuition of mass function as measuring the
effect of time which is a solution to the scalar wave equation in this
setting. We mention also the fact that in temporal gauge the effect
or variations of time $t$ are visible by all observers in the
spacetime $(g(\nabla t,\nabla t)=g^{-1}(dt,dt)<0)$, and our collapsing process
includes the assignment of a coordinates system which assures: the
compactification of the spacetime, that time effect ceases to be
visible at later times, and time collapses to zero.
\section{Compactification of Spacetime}
We assume independently of the matter contents that a global in time
spacetime $(V,g)$ is established by solving the Einstein-field
equations in temporal gauge, with $g$ of the following form in
coordinates $(x^\alpha)$:
\begin{equation}
    g=-\frac{|\overline{g}|}{\Omega^2}(dx^0)^2+\overline{g}_{ij}dx^idx^j.
\end{equation} $\Omega$ is a $x^0$-independent function,
\begin{equation}
    \Omega^2 g=-|\overline{g}| (dx^0)^2 +\Omega^2\overline{g}_{ij}
    dx^idx^j.
\end{equation}
We assume that $\Omega>0$ is a bounded quantity tending to $0$ as
$r:=\sqrt{\displaystyle\sum_{i=1}^n(x^i)^2}$ tends
 to $+\infty$,
since it represents some initial scalar density in the spacetime.
Furthermore, the gradient $\nabla \Omega$ does not vanish and in
particular we admit that $\frac{\partial\Omega}{\partial x^1}$ does
not vanish.
\subsection{New gauge or observer in an accelerated reference frame}
 We proceed to a
rescaling of the spacetime defining $\Phi$:
\begin{equation}\label{22}
  \begin{array}{ll}
   \widehat{ \Phi} : V\rightarrow \widehat{V}=]0;1]\times ]0;
   \Omega_{\max}]\times (]-\frac{\pi}{2};+\frac{\pi}{2}[)^{n-1}  & \hbox{} \\
     x=(x^\alpha) \mapsto
  \omega:=\left(e^{-\Omega\;x^0},\Omega,\arctan(\Omega\;x^a)\right):=(\omega^\alpha). & \hbox{}
  \end{array}
\end{equation}One has the following expressions for $i,s=1,...,n; a,b=2,...,n$:
\begin{eqnarray*}\label{}
    \frac{\partial \omega^0}{\partial x^0}=-\Omega\;\omega^0
  ,\; \frac{\partial \omega^s}{\partial x^0}=0,
  \frac{\partial \omega^0}{\partial x^i}=-x^0\omega^0\frac{\partial\Omega}{\partial x^i}
  , \;  \frac{\partial \omega^1}{\partial x^i}=\frac{\partial\Omega}{\partial x^i},&&\\
  \frac{\partial \omega^a}{\partial x^1}=\frac{x^a\frac{\partial\Omega}{\partial x^1}}{1+(\Omega\;x^a)^2},\;\frac{\partial \omega^a}{\partial x^b}= \frac{x^a\frac{\partial\Omega}{\partial x^b}+\Omega\delta_b^a}{1+(\Omega\;x^a)^2}.&&
\end{eqnarray*}These expressions induce that the change of variables $\Phi:\;x^\alpha \rightarrow \omega^\delta$ observes the following jacobian:
\begin{equation}\label{21}
   J\equiv \left|\frac{D(\omega)}{D(x)}\right|=-\omega^0\Omega^n\frac{\partial\omega^1}{\partial
    x^1}\prod_{a=2}^n\left(\frac{1}{1+\tan^2\omega^a}\right),
\end{equation}

The inverse of the jacobian matrix has the following components:
\begin{eqnarray*}
    \frac{\partial x^0}{\partial \omega^0}=-\frac{1}{\Omega\;\omega^0}
  ,\;
  \frac{\partial x^0}{\omega^1}=-\frac{x^0}{\Omega},\;\frac{\partial x^0}{\omega^a}=0,
     &&\\
  \frac{\partial x^1}{\partial\omega^0}=0,\;\frac{\partial x^1}{\partial \omega^1}=\frac{\Omega+\sum_{a=2}^{n}x^a\frac{\partial\Omega}{\partial x^a}}{\Omega\frac{\partial\Omega}{\partial x^1}},\;
   \frac{\partial x^1}{\partial \omega^a }=-\frac{\frac{\partial\Omega}{\partial x^a}(1+\tan^2\omega^a)}{\Omega\frac{\partial\Omega}{\partial x^1}},&&\\
 \frac{\partial x^a}{\partial \omega^0}=0,\; \frac{\partial x^a}{\partial\omega^1}=- \frac{x^a}{\Omega},\; \frac{\partial x^a}{\partial\omega^a}=\frac{1+\tan^2\omega^a}{\Omega},\;\frac{\partial x^a}{\partial\omega^b}=0,\;a\neq b.&&
\end{eqnarray*}

 {\remark{}\label{r2} The manifold $\widehat{V}$
has two specific boundary hypersurfaces $i^+:\;\omega^0=0$ and
$\mathcal{I}:\;\omega^1=0$, correspondingly codimension $2$ corners
coexist.}
\subsection{ b-Lorentzian metric -- Scattering Lorentzian metric -- Scattering dual metric}
The change of coordinates from $V$ to a manifold with corners
$\widehat{V}$ induces the necessity to using the $b-geometry$ on
manifolds with corners \cite{96}, $"b"$ referring to boundary. We
only mention that the manifold with corners $\widehat{V}$ infers a
Lie algebra $\mathcal{V}(\widehat{V})$ of vector fields on
$\widehat{V}$ tangent to the boundary called $b-vector\; fields$ as
usual. $\mathcal{V}(\widehat{V})$ is the space of sections of a
natural vector bundle $^bT\widehat{V}$ over $\widehat{V}$, the
$b-tangent\; bundle$, which over the interior of $\widehat{V}$ is
naturally identified with $T\widehat{V}$. $^bT\widehat{V}$ is
spanned near the two-boundary hypersurfaces $\omega^0=0$ and
$\omega^1=0$ by
$\omega^0\frac{\partial}{\partial\omega^0},\omega^1\frac{\partial}{\partial\omega^1},\frac{\partial}{\partial\omega^a}$,
and $\mathcal{V}(\widehat{V})$ is spanned over
$\mathcal{C}^\infty(\widehat{V})$ by these vector fields. The bundle
dual of $^bT\widehat{V}$ denoted $^bT^\ast\widehat{V}$ and called
$b-cotangent \;bundle$, is spanned locally near the two-boundary
hypersurfaces $\omega^0=0$ and $\omega^1=0$ by
$\frac{d\omega^0}{\omega^0},\frac{d\omega^1}{\omega^1},d\omega^a$.\\
    Relatively to this dual basis, the unphysical metric $\Omega^2g\equiv h$ has the form
   \begin{equation*} h=\widehat{h}_{00}\left(\frac{d\omega^0}{\omega^0}\right)^2+\widehat{h}_{01}
    \left(\frac{d\omega^0}{\omega^0}\otimes \frac{d\omega^1}{\omega^1}+
    \frac{d\omega^1}{\omega^1}\otimes
    \frac{d\omega^0}{\omega^0}\right)+\widehat{h}_{11}\left(\frac{d\omega^1}{\omega^1}\right)^2+
\end{equation*}
\begin{equation}\label{20}
\widehat{h}_{1a}
    \left(\frac{d\omega^1}{\omega^1}\otimes d\omega^a+
    d\omega^a\otimes
    \frac{d\omega^1}{\omega^1}\right)+\widehat{h}_{ab}d\omega^a\otimes
    d\omega^b,
\end{equation}which we refer to as the b-Lorentzian metric.
\begin{equation}
  \widehat{h}_{00}=- \frac{|\overline{g}|}{\Omega^2},\;\widehat{h}_{01}=|\overline{g}|\frac{\ln\omega^0}{\Omega^2};
\end{equation}
\begin{equation}
 \begin{split}
     \widehat{h}_{11}=& -|\overline{g}|\left(\frac{\ln\omega^0}{\Omega}\right)^2+\Omega^4g_{11}
  \left(\frac{\Omega+\sum_{a=2}^{n}x^a\frac{\partial\Omega}{\partial x^a}}{\Omega\frac{\partial\Omega}{\partial x^1}}\right)^2- \\
      & 2\Omega^2\left(\frac{\Omega+\sum_{a=2}^{n}x^a\frac{\partial\Omega}{\partial x^a}}{\Omega\frac{\partial\Omega}{\partial x^1}}\right)\sum_{c=2}^{n} g_{1c}\tan\omega^c+g_{ab}\tan\omega^a\tan\omega^b;
 \end{split}
\end{equation}
\begin{equation}
\begin{split}
    \widehat{h}_{1a}= & -\Omega^3 g_{11}\left(\frac{\Omega+\sum_{a=2}^{n}x^a\frac{\partial\Omega}{\partial x^a}}{\Omega\frac{\partial\Omega}{\partial x^1}}\right)
  \left(\frac{\frac{\partial\Omega}{\partial x^a}(1+\tan^2\omega^a)}{\Omega\frac{\partial\Omega}{\partial x^1}}\right)+ \\
     & \Omega^2g_{1a}\left(\frac{\Omega+\sum_{c=2}^{n}x^c\frac{\partial\Omega}{\partial x^c}}{\Omega\frac{\partial\Omega}{\partial x^1}}\right)(1+\tan^2\omega^a)+ \\
     &
  \Omega\sum_{c=2}^{n}g_{1c}\tan\omega^c\left(\frac{\frac{\partial\Omega}{\partial x^a}(1+\tan^2\omega^a)}{\Omega\frac{\partial\Omega}{\partial x^1}}\right)-
   \sum_{c=2}^{n}g_{ac}(1+\tan^2\omega^a)\tan\omega^c;
\end{split}
\end{equation}
\begin{equation}
  \begin{split}
    \widehat{h}_{ab}=& g_{11}\frac{(1+\tan^2\omega^a)(1+\tan^2\omega^b)}{\left(\frac{\partial\Omega}{\partial x^1}\right)^2}
  \left(\frac{\partial\Omega}{\partial x^a}\right)\left(\frac{\partial\Omega}{\partial x^b}\right)+ \\
       & 2g_{1a}\frac{(1+\tan^2\omega^a)(1+\tan^2\omega^b)}{\left(\frac{\partial\Omega}{\partial x^1}\right)}\left(\frac{\partial\Omega}{\partial x^b}\right)+g_{ab}(1+\tan^2\omega^a)(1+\tan^2\omega^b).
  \end{split}
\end{equation}
To this
b-Lorentzian metric corresponds the volume density $dh$, a
non-vanishing b-density and one can better consider b-density bundles
in general which are locally spanned by $\frac{d\omega^0}{\omega^0}
\frac{d\omega^1}{\omega^1}d\omega^2...d\omega^n$.\\
The scattering Lorentzian metric in this context in turn corresponds to the writing of $g$ and its dual with respect to the bases
$\left(\frac{d\omega^0}{\Omega\omega^0},\;\frac{d\Omega}{\Omega^2},\;\frac{d\omega^a}{\Omega}\right)$,  $(\omega^0\Omega\partial_{\omega^0},\;\Omega^2\partial_\Omega ,\;\Omega\omega^a)$.
\\The scattering dual metric reads:
\begin{equation}\label{a1}
\begin{split}
   g^{-1}= & \left(g^{00}+\frac{(\ln\omega^0)^2}{\Omega^4}g^{ij}\frac{\partial\Omega}{\partial x^i}
   \frac{\partial\Omega}{\partial x^j}\right)(\Omega\omega^0\partial_{\omega^0})^2+2\frac{\ln\omega^0}{\Omega^4}g^{ij}\frac{\partial\Omega}{\partial x^i}
   \frac{\partial\Omega}{\partial x^j}\Omega\omega^0\partial_{\omega^0}\otimes\Omega^2\partial_{\Omega}+ \\
     & 2\frac{\ln\omega^0}{\Omega^3}\frac{\partial\Omega}{\partial x^i}\left(\frac{\partial\Omega}{\partial x^j}\frac{\tan\omega^a}{\Omega}+\Omega\delta_j^a\right)g^{ij}
     \Omega\omega^0\partial_{\omega^0}\otimes\Omega\partial_{\omega^a}+\frac{1}{\Omega^4}g^{ij}\frac{\partial\Omega}{\partial x^i}
   \frac{\partial\Omega}{\partial x^j}(\Omega^2\partial_{\Omega})^2+\\
     & 2\frac{\partial\Omega}{\partial x^i}\frac{\left(\frac{\partial\Omega}{\partial x^j}\frac{\tan\omega^a}{\Omega}+\Omega\delta_j^a\right)}{\Omega^3(1+\tan^2\omega^a)}g^{ij}
     \Omega^2\partial_{\Omega}\otimes\Omega\partial_{\omega^a}+\\
     & \frac{1}{\Omega^2}\frac{(\frac{\partial\Omega}{\partial x^i}\frac{\tan\omega^a}{\Omega}+\Omega\delta_i^a)
     (\frac{\partial\Omega}{\partial x^j}\frac{\tan\omega^b}{\Omega}+\Omega\delta_j^b)g^{ij}}{(1+\tan^2\omega^a)(1+\tan^2\omega^b)}
     \Omega\partial_{\omega^a}\otimes \Omega\partial_{\omega^b}.
\end{split}
  \end{equation}
{\remark{} The
b-Lorentzian metric (\ref{20}) together with its dual
and also the scattering Lorentzian metric and its dual are very useful in studying asymptotic
questions to the two boundary hypersurfaces $\omega^0=0$, and
$\omega^1=0$, notably the analysis of the structure of the
null-geodesics flow, that is the flow of the Hamilton function
within the characteristic set \cite{43}. In the gauge here it appears from (\ref{a1}) that the regularity of infinity depends on the estimate of $g^{ij}\frac{\partial\Omega}{\partial x^i}
   \frac{\partial\Omega}{\partial x^j}$.\\
 We
concentrate in this article on a lower bound inequality for the "total mass"
assuming that the considered system collapses to a black hole.   }\\
For what follows, some definitions are necessary. We signal that
some terminologies used below might not coincide with the usual one
in the literature.
\section{Mass function and derived geometric analysis}
\subsection{Motivations}
 The conformal geometry here in order to fit to collapse
requires that relative time effect passes from visible at early
times to invisible at later ones. That is
$\frac{d\omega^0}{\omega^0}$ passes from the interior of the dual
light cone to its exterior.
\\The dual of the conformal metric $h\equiv\Omega^2
g=-|\overline{g}|(dx^0)^2 +\Omega^2\overline{g}_{ij}dx^idx^j$ in
coordinates $( x^\alpha)$ reads:
\begin{equation}
    G\equiv h^{-1}=-\frac{1}{|\overline{g}|}\left(\frac{\partial}{\partial
    x^0}\right)^2+\Omega^{-2}\overline{g}^{ij}\frac{\partial}{\partial x^i}\frac{\partial}{\partial
    x^j}.
\end{equation}In coordinates $(\omega^\alpha)$,
 components of tensors are decorated with a $hat$  $" \;\widehat{}\;"$,
one has in particular:
\begin{equation}
   G(d\omega^0,d\omega^0)= \widehat{G}^{00}=(\omega^0)^2\Omega^2\left(-\frac{\Omega^2}{|\overline{g}|}+
    \frac{(\ln\omega^0)^2}{\Omega^4}\overline{g}^{ij}\frac{\partial\Omega}{\partial x^i}\frac{\partial
    \Omega}{\partial
    x^j}\right).
\end{equation}On the hypersurface $\omega^0=1$, $\widehat{G}^{00}=
-\frac{\Omega^4}{|\overline{g}|}_{|\omega^0=1}<0$, provided $\frac{1}{\Omega^2}\overline{g}^{ij}\frac{\partial\Omega}{\partial x^i}\frac{\partial
    \Omega}{\partial
    x^j}_{|\omega^0=1} \nrightarrow +\infty$; the
covariant vector $\frac{d\omega^0}{\Omega\omega^0}$ is timelike near
$\omega^0=1$. Now, conceptually it is desirable that  this
covariant vector $\frac{d\omega^0}{\Omega\omega^0}$ is spacelike near the hypersurface $\omega^0=0$, and that
the equation  $G(\frac{d\omega^0}{\Omega\omega^0},\frac{d\omega^0}{\Omega\omega^0})=0$ defines a smooth hypersurface. This induces certainly some estimates for the term  $\frac{(\ln\omega^0)^2}{\Omega^4}\overline{g}^{ij}\frac{\partial\Omega}{\partial x^i}\frac{\partial
    \Omega}{\partial
    x^j}$ that would eventually be elucidated according to the desired end structure at infinity. Let's recall that in the orthogonal splitting of the metric $g= -\alpha^2(dt)^2+\overline{g}$, the lapse $\alpha=\frac{\sqrt{|\overline{g}|}}{\Omega}$ can be chosen as bounded (see O. M$\ddot{u}ller$ ,  M. Sanchez, 2009, and references in) so that the behavior of
$\frac{\widehat{G}^{00}}{(\Omega\omega^0)^2}$ at the boundaries hypersurfaces is given by the one of
$\frac{(\ln\omega^0)^2}{\Omega^4}\overline{g}^{ij}\frac{\partial\Omega}{\partial x^i}\frac{\partial
    \Omega}{\partial
    x^j}$. On the other
hand
$\widehat{G}^{11}=(\omega^1)^{-2}\overline{g}^{ij}\frac{\partial\omega^1}{\partial
x^i}\frac{\partial\omega^1}{\partial x^j}$ and since
$\nabla\Omega\neq (0)$ one has $\widehat{G}^{11}>0$ and then
$\frac{d\omega^1}{\omega^1}$ is spacelike for $\omega^1\neq 0$.
Depending on the chosen end structure, $\omega^1=0$ might be a null
hypersurface. \\According to what precedes and assuming some smoothness assumptions, there exits a hypersurface beyond which $\frac{d\omega^0}{\Omega\omega^0}$
ceases to be timelike. It seems that time ceases to exist on such
hypersurface. As indicated by the theory of Relativity due to Sir
Albert Einstein, the speed of light is the same everywhere whether
one is under the influence of very high or low gravity. So it is
reasonable that something decreases as one approaches the speed of
light, and it is the time or its effect. The same situation is valid
for a Black hole since light can not escape from it so the time
nearby it is hardly visible and felt. It is difficult to notice the
change in the Black hole's surroundings. We then resume that the
manifold with corners $\widehat{V}$ is split in two (exterior and
interior) regions as $\frac{d\omega^0}{\Omega\omega^0}$ moves from the
interior of the dual characteristic cone to its exterior. We
therefore adopt the following definitions. {\definition {A system of
coordinates $(\omega^\alpha)$ of the manifold with corners
$(\widehat{V},h)$ is called gravitational
collapse-adapted if one of the coordinates $\omega^\alpha$, say $\omega^0$, satisfies:\\
- the vector field $\frac{\partial}{\partial\omega^0}$ is timelike
 at least beyond some hypersurface (or asymptotically timelike,
 i.e. timelike in the neighborhood of the final state $\omega^0 =0$),\\
- the dual vector $d\omega^0$ passes from the timelike character in
some exterior region to the spacelike character in the neighborhood
of the final state.}}
\subsection{Heuristics}
 To motivate more our present understanding of
the event horizon, let's recall that the local energy density in a
spacetime equipped with a system of coordinates $(x^\alpha)$ is
defined as $\mu=\frac{1}{8\pi}E^{00}$, where
 $E$ is the Einstein tensor. In a spacetime solution of the
  Einstein's equations with matter, this local energy
 equals $T^{00}$. $T$ is the stress energy-momentum
  tensor of matter. One can then write:
 \begin{equation*}
   \frac{R^{00}-8\pi T^{00}}{R} = \frac{1}{2}g^{00}.
 \end{equation*}
 In relation with the gauge $(\omega^\alpha)$, one has
 \begin{equation}\label{a2}
 \begin{split}
    G\left(\frac{d\omega^0}{\Omega\omega^0}, \frac {d\omega^0}{\Omega\omega^0}\right) & =
    g^{-1}(dx^0,dx^0)+
    \frac{(\ln\omega^0)^2}{\Omega^4}\overline{g}^{ij}\frac{\partial\Omega}{\partial x^i}\frac{\partial
    \Omega}{\partial
    x^j} \\
      & =\frac{2}{R}(R^{00}-8\pi T^{00})+
     \left(\frac{x^0}{\Omega}\right)^2\overline{g}^{ij}\frac{\partial\Omega}{\partial x^i}\frac{\partial
    \Omega}{\partial
    x^j}  \\
      & = h\left(\frac{\widehat{\nabla}\omega^0}{\Omega\omega^0}, \frac {\widehat{\nabla}\omega^0}{\Omega\omega^0}\right)).
 \end{split}
    \end{equation}These expressions are useful in various respects. The first equality could be used in vacuum.
 They suggest that the region where time effects are
 hidden
  is described locally by $ G\left(\frac{d\omega^0}{\Omega\omega^0}, \frac {d\omega^0}{\Omega\omega^0}\right)>0$
   and corresponds to some interior region
   whereas the region $ G\left(\frac{d\omega^0}{\Omega\omega^0}, \frac {d\omega^0}{\Omega\omega^0}\right)<0$ is an
 exterior region where time remains visible. One can also claim that in vacuum every gravitational collapse is dictated
  by the time effect. In the case where $G\left(\frac{d\omega^0}{\Omega\omega^0}, \frac {d\omega^0}{\Omega\omega^0}\right)$
  tends to $+\infty$ as $\omega^0$ tends to zero,
 we call $\omega^0=0$ a mass singularity. Finally the total time effect
 is expected to be responsible to the final state of the system.\\
  The hypersurface $G\left(\frac{d\omega^0}{\Omega\omega^0}, \frac {d\omega^0}{\Omega\omega^0}\right)=0$ if smooth defines consequently an
  (illusory) event horizon which might
 comprise a part of the blackhole's event horizon and
  a part equals to the apparent horizon.
  \subsection{Geometric objects and properties}
    {\definition The mass function at a point $(\omega^\alpha)$
in the manifold with corners $(\widehat{V},h\equiv \Omega^2 g)$
is the quantity $m(\omega^0,\omega^1,...,\omega^n)\equiv
G(\frac{d\omega^0}{\Omega\omega^0},
\frac{d\omega^0}{\Omega\omega^0})(\omega^0,\omega^1,...,\omega^n)$.\\
We remark that this definition of the mass is valid up to the
boundary of the considered manifold since
$\frac{d\omega^0}{\omega^0}$ is non-singular (and non-trivial) as a
b-covector at $\omega^0=0$, and the same holds for
$\frac{d\omega^1}{\omega^1}$ at $\omega^1=0$. This mass-function
measures some extra local energy density since it may exist even in
vacuum. At every point $(\omega^\alpha)$ it encodes the rate at
which the dual vector $d\omega^0$ passes from the interior of the
dual characteristic cone at $(\omega^\alpha)$ to its exterior.}
{\notation{} In the rest of the paper we adopt the notation
$d\varpi=d\omega^2...d\omega^n$.}
 {\definition{} The (squared) mass of a
region $D$ in the the manifold with corners $(\widehat{V},h)$ is the
quantity
\\$\int_D
m(\omega^0,\omega^1,...,\omega^n)d\omega^0d\omega^1d\varpi$.
 This quantity may be negative meaning that its squared root is a complex number.}
  {\definition{} The (squared) total mass of the spacetime is the quantity\\
$\int_{\widehat{V}}
m(\omega^0,\omega^1,...,\omega^n)d\omega^0d\omega^1d\varpi$. }
{\definition{} The (illusory) Event horizon in the spacetime $(V,g)$ is the
hypersurface $\mathcal{H}$ of points $(\omega^\alpha)$ in the
manifold with corners $(\widehat{V},h)$ which satisfy the equation
\begin{equation}
m(\omega^0,\omega^1,...,\omega^n)=0.
\end{equation}
 }
 Provided a monotonicity argument for $m$ is established as described in the sequel of the paper,
  the existence of the (illusory) event horizon ( assumed smooth for simplicity ) induces the following definitions.
{\definition{} The interior region in the spacetime $(V,g)$ denoted
$Int (V)$ is:
\begin{equation}
Int (V)=\{\omega\in
\widehat{V}/m(\omega^0,\omega^1,...,\omega^n)\geq 0\}.
\end{equation}
 }
 {\definition{} The exterior region in the spacetime $(V,g)$ denoted
$Ext (V)$ is:
\begin{equation}
Ext (V)=\{\omega\in \widehat{V}/m(\omega^0,\omega^1,...,\omega^n)<
0\}.
\end{equation}
 }
{\remark{} The wave fronts (WF) for a system of partial differential
equations are submanifolds $\{f=cste\}$ of the spacetime where the
scalar function $f$ is solution of the eikonal equation (\cite{11},
P.\; 270). A wave front is generated by the bicharacteristics or rays (the characteristics of the eikonal equation).\\ In this setting (i.e. Temporal gauge), the principal part
in the considered evolution system is the operator $\square
\partial_0$, the characteristic determinant w.r.t. the coordinates
$(x^\alpha)$ is the hyperbolic polynomial
\begin{equation*}
    P(x,\zeta)=\left(g^{00}(\zeta_0)^2+g^{ij}\zeta_i\zeta_j\right)\zeta_0,
\end{equation*}and
the corresponding eikonal equation is
\begin{equation*}
  \frac{\partial f}{\partial x^0}  G^{\alpha\beta}\frac{\partial f}{\partial x^\alpha}\frac{\partial
    f}{\partial x^\beta}=0.
\end{equation*}Furthermore the change of coordinates $(x\rightarrow \omega)$ induces that
\begin{equation*}
    \square\partial_0 \simeq -\Omega\omega^0 \widehat{\square}
    \widehat{\partial}_0,
\end{equation*}correspondingly
\begin{equation*}
    P(\omega,\zeta)=-\Omega\omega^0 \widehat{G}^{\alpha\beta}\zeta_\alpha\zeta_\beta\zeta_0.
\end{equation*}Two types of wave fronts are then on concerned here:\\
- $f=cste$ s.t. $\frac{\partial f}{\partial\omega^0}=0$,\\
- $f=cste$ s.t. $ \widehat{G}^{\alpha\beta}\frac{\partial
f}{\partial\omega^\alpha}\frac{\partial
    f}{\partial\omega^\beta}=0$.\\
   Given a point $\omega\in \widehat{V}$, the tangents to the rays issued from $\omega$ generate a cone in the tangent space to the spacetime at $\omega$, called the wave cone, dual to the characteristic cone. This wave cone is the envelope of the hyperplanes whose normals belong to the characteristc cone.\\ We denote by $WF (\omega)$ the wave front at the point $\omega$.
    $WF (\omega) $ is tangent to the wave cone at $\omega$ along the
    direction of the bicharacteristic or ray. }
  {\definition{} The Apparent event horizon of the spacetime $(V,g)$
is the surface denoted $\mathcal{H}^-$ and defined by:
\begin{equation}
\mathcal{H}^-=\{\omega\in \mathcal{H}/ \;WF (\omega)\cap Ext (V)\neq
\emptyset\}.
\end{equation}}
 {\definition{} The actual event horizon of the spacetime
is the surface denoted $\mathcal{H}^+$ and defined by:
\begin{equation}
\mathcal{H}^+=\{\omega\in \mathcal{H}/ \; WF (\omega)\subset Int
(V)\}.
\end{equation}}It is obvious that $\mathcal{H}=\mathcal{H}^+\cup
\mathcal{H}^-$. {\definition{} The Black hole region of the
spacetime $V$ denoted $\mathcal{\textbf{B}}$ is defined as:
\begin{equation}
 \mathcal{\textbf{B}}=\{\omega\in Int (V)/ WF (\omega)\subset Int (V)\}.
\end{equation}}
{\definition{} The Blackhole's event horizon is the boundary of the
Black hole region $\mathcal{\textbf{B}}$. We denote it $\partial
\mathcal{\textbf{B}}$. }
 {\definition{} The White hole region of the
spacetime denoted $\mathcal{\textbf{W}}$ is:
\begin{equation}
\mathcal{\textbf{W}}=Int (V)\setminus \mathcal{\textbf{B}}.
\end{equation}Its boundary is denoted $\partial \mathcal{\textbf{W}}$}
  {\remark{ Under the decreasing property for $m$ w.r.t. $\omega^0$:\\
 a- the illusory event horizon observes obviously
 the following disjoint
union:
\begin{equation}
\mathcal{H}=\{\omega=(\omega^0,\omega^i)\in \mathcal{H}/
\frac{\partial m}{\partial \omega^0}(\omega)< 0\}\cup
\{\omega=(\omega^0,\omega^i)\in \mathcal{H}/ \frac{\partial
m}{\partial \omega^0}(\omega)=0\};
\end{equation}
b- the subset $\{m(\omega^0,\omega^i)=0,\;\frac{\partial
m}{\partial\omega^0}(\omega^0,\omega^i)<0\}$ of $\mathcal{H}$ has an equation of the form $ \omega^0= X(\omega^i)$, by the implicit function theorem.\\
 c- the condition for the event horizon $\mathcal{H}$
  to be a wave front is that $m$ satisfies
  the eikonal equation, i.e. $\frac{\partial
m}{\partial\omega^0}\;\widehat{G}^{\lambda\delta}\frac{\partial
m}{\partial \omega^\lambda}\frac{\partial m}{\partial
\omega^\delta}=0$. This equation reduces if $\mathcal{H}$ has an equation $\omega^0= X(\omega^i)$ to
\begin{equation}
    \left[\frac{\partial m}{\partial\omega^0}\right]^3\left(-2\left[\widehat{G}^{0i}\right]\frac{\partial
    X(\omega^i)}{\partial\omega^i}+\left[\widehat{G}^{ij}\right]\frac{\partial
    X(\omega^i)}{\partial\omega^i}\frac{\partial
    X(\omega^i)}{\partial\omega^j}\right)=0;
\end{equation}
d- It is reasonable to assume that
any blackhole's event horizon admits an equation of the form $
\omega^0-\widetilde{X}(\omega^i)=0$ since it must form a "hoop".
Intuitively the singularity here is pushed off to infinity, indeed, $G(d\omega^0, d\omega^0)=0$ if and only if
$\omega^0=0$ or $\Omega=0$ or $m=0$. }}
{\proposition{a-The (illusory) event horizon $\mathcal{H}$ is not necessarily a
null surface. \\
b- At each point $\omega$ of the actual event horizon
$\mathcal{H}^+$, one has $n$-dimensional wave fronts
$\{\omega^0=cste\}$.\\
c- The apparent horizon $\mathcal{H}^-$ is a wavefront, it has always
a null
character.\\
d- If the (illusory) event horizon $\mathcal{H}$ admits an equation
of the form $\omega^0=cste$, then it is a null hypersurface.
}}\\
\textbf{Proof}\\Obvious since the (illusory) event horizon is of equation
$m\equiv\frac{\widehat{G}^{00}}{(\omega^0)^2\Omega^2}=0$ and the eikonal
equation reads $\frac{\partial
f}{\partial\omega^0}\widehat{G}^{\mu\nu}\frac{\partial
f}{\partial\omega^\mu}\frac{\partial f}{\partial\omega^\nu}=0$.
Furthermore, $\mathcal{H}^-=\{m=0,\;\frac{\partial
m}{\partial\omega^0}=0\}$.
\subsection{Monotonicity property for the mass function}
 {\theorem {}\label{th2}
  The mass function $m$ is decreasing relatively to the time $\omega^0$ if and only if
the extrinsic curvature $K$ of the spacetime satisfies the
inequality
\begin{equation}\label{16}
trK+\frac{2|\overline{g}|\sqrt{|\overline{g}|}\ln\omega^0}{\Omega^4}g^{ij}\frac{\partial\Omega}{\partial x^i}\frac{\partial\Omega}{\partial x^j}-\frac{2|\overline{g}|^2(\ln\omega^0)^2}{\Omega^6}K^{ij}\frac{\partial\Omega}{\partial x^i}\frac{\partial\Omega}{\partial x^j}\leq 0.
\end{equation}}
 \textbf{Proof}\\
 The proof of this theorem is based merely on straightforward
computations. We recall for this purpose that $m$ is given by:
\begin{equation*}
    m(\omega^\alpha)=\left(-
    \frac{\Omega^2}{|\overline{g}|}+
    \frac{(\ln\omega^0)^2}{\Omega^4}\overline{g}^{ij}\frac{\partial\Omega}{\partial x^i}\frac{\partial
    \Omega}{\partial
    x^j}\right).
\end{equation*} In order to address the expression of $\frac{\partial m}{\partial\omega^0}$, we
collect the following useful expressions:
\begin{equation}
    \frac{\partial
    x^0}{\partial\omega^0}=-\frac{1}{\Omega\omega^0},\;\frac{\partial|\overline{g}|}{\partial
    \omega^0}=-\frac{1}{\Omega\omega^0}\frac{\partial|\overline{g}|}{\partial
    x^0},\;
    \frac{\partial g^{ij}}{\partial
    \omega^0}=-\frac{1}{\Omega\omega^0}\frac{\partial g^{ij}}{\partial
    x^0}.
\end{equation}Under the temporal gauge condition, one has
\begin{equation*}
    \frac{\partial \overline{g}_{ij}}{\partial
    x^0}=-2\frac{\sqrt{|\overline{g}|}}{\Omega}K_{ij},\;
    \frac{\partial \overline{g}^{ij}}{\partial
    x^0}=2\frac{\sqrt{|\overline{g}|}}{\Omega}K^{ij},\;
    \frac{\partial|\overline{g}|}{\partial
    x^0}=\frac{1}{2}\overline{g}^{ij}\frac{\partial \overline{g}_{ij}}{\partial
    x^0};
\end{equation*}where $K=(K_{ij})$ is the extrinsic curvature of the
hypersurfaces $x^0=t$. It follows that:
\begin{equation}
\frac{\partial|\overline{g}|}{\partial
    \omega^0}=\frac{1}{\Omega^2\omega^0}\sqrt{|\overline{g}|}trK,\;
    \frac{\partial \overline{g}^{ij}}{\partial\omega^0}=
    -\frac{2}{\Omega^2\omega^0}\sqrt{|\overline{g}|}K^{ij}.
\end{equation}Combining all these expressions one has:
\begin{equation}\label{5}
    \frac{\partial
    m}{\partial\omega^0}=\frac{1}{\omega^0}\left(\frac{\sqrt{|\overline{g}|}}{|\overline{g}|^2}trK+
    \frac{2\ln\omega^0}{\Omega^4}g^{ij}\frac{\partial\Omega}{\partial x^i}\frac{\partial\Omega}{\partial x^j}-\frac{2\sqrt{|\overline{g}|}(\ln\omega^0)^2}{\Omega^6}K^{ij}\frac{\partial\Omega}{\partial x^i}\frac{\partial\Omega}{\partial x^j}\right).
\end{equation}The conclusion is then obvious according to this
expression.\\
{\corollary\label{co1}
  For the mass function to be
strictly decreasing w.r.t. to the time $\omega^0$, it suffices that
 the extrinsic curvature $K$ of the spacetime satisfies the
inequalities
\begin{equation}
trK < 0,\; K^{ij}\frac{\partial\Omega}{\partial x^i}\frac{\partial
   \Omega}{\partial
    x^j}\geq
    0.
\end{equation}
  }
  \textbf{Proof}\\
  It is obvious according to the expression (\ref{5}) of $\frac{\partial
  m}{\partial\omega^0}$ above.
  {\corollary\label{co11}
  For the mass function to be
strictly decreasing w.r.t. to the time $\omega^0$, it suffices that
 the extrinsic curvature $K$ of the spacetime satisfies the
inequalities
\begin{equation}
(lapse)^2\left(g^{ij}\frac{\partial\Omega}{\partial x^i}\frac{\partial
   \Omega}{\partial
    x^j}\right)^2+2trK \;K^{ij}\frac{\partial\Omega}{\partial x^i}\frac{\partial
    \Omega}{\partial
    x^j}< 0,\; K^{ij}\frac{\partial\Omega}{\partial x^i}\frac{\partial
    \Omega}{\partial
    x^j}>
    0.
\end{equation}
  }
  \textbf{Proof}\\
  The proof is obvious according to the expression (\ref{5}) of $\frac{\partial
  m}{\partial\omega^0}$ above.

  {\remark{In the rest of this work, it would be interesting to see how $\frac{\partial m}{\partial \omega^0}$ is related to
  $\frac{dm}{ds}$ where $s$ is the black hole's horizon perturbation parameter.}}
  \subsection{Monotonicity property: Another view}
  In the following lines, we are interested in another possible geometrical
analysis and physical interpretation of the monotonicity property
(\ref{16}) above, i.e.:
\begin{equation}\label{37}
    \frac{\partial m}{\partial\omega^0}\leq 0;
\end{equation}
  and eventually the condition (\ref{27})
used in the sequel of this work, i.e.:
\begin{equation}\label{38}
    \widehat{h}^{0i}\frac{\partial
    m}{\partial\omega^i}>0,\;\hbox{on}\;
    \mathcal{H}^+.
\end{equation}
We address the following computations where covariant derivatives
are carried according to the unphysical metric $h$. Given then
$\widehat{\nabla}$ the riemannian connection of $(\widehat{V},h)$,
one has:
\begin{equation}
\frac{\partial m}{\partial\omega^0}=\frac{\partial }{\partial\omega^0}
  \left(h\left(\frac{\widehat{\nabla}\omega^0}{\Omega\omega^0},\frac{\widehat{\nabla}\omega^0}{\Omega\omega^0}\right)\right)
=2h\left(\widehat{\nabla}_{\frac{\partial}{\partial\omega^0}}^{\frac{\widehat{\nabla}\omega^0}{\Omega\omega^0}},
\frac{\widehat{\nabla}\omega^0}{\Omega\omega^0}\right)
\end{equation}From this expression, it results that $\frac{\partial
m}{\partial\omega^0}<0$, if and only if
\begin{equation}\label{40}
h\left(\widehat{\nabla}_{\frac{\partial}{\partial\omega^0}}^{\frac{\widehat{\nabla}\omega^0}{\Omega\omega^0}},
\frac{\widehat{\nabla}\omega^0}{\Omega\omega^0}\right) <0.
\end{equation}On the other hand and similarly with the above computations, one has:
\begin{equation}\label{41}
\frac{\partial m}{\partial\omega^i}=2h\left(\widehat{\nabla}_{\frac{\partial}{\partial\omega^i}}^{\frac{\widehat{\nabla}\omega^0}{\Omega\omega^0}},
\frac{\widehat{\nabla}\omega^0}{\Omega\omega^0}\right),
\end{equation}this implies that:
\begin{equation}\label{42}
\widehat{h}^{0i}\frac{\partial m}{\partial\omega^i}=2
h\left(\widehat{\nabla}_{\widehat{h}^{0i}\frac{\partial}{\partial\omega^i}}^{\frac{\widehat{\nabla}\omega^0}{\Omega\omega^0}},
\frac{\widehat{\nabla}\omega^0}{\Omega\omega^0}\right).
\end{equation}
{\proposition{}\label{p3}i- The condition $\frac{\partial
m}{\partial\omega^0}<0$ is equivalent to
\begin{equation}\label{39}
   h\left(\widehat{\nabla}_{\frac{\partial}{\partial\omega^0}}^{\frac{\widehat{\nabla}\omega^0}{\Omega\omega^0}},
\frac{\widehat{\nabla}\omega^0}{\Omega\omega^0}\right) <0.
\end{equation}
ii- The condition $\widehat{h}^{0i}\frac{\partial
m}{\partial\omega^i}>0$ is equivalent to \begin{equation}\label{43}
h\left(\widehat{\nabla}_{\widehat{h}^{0i}\frac{\partial}{\partial\omega^i}}^{\frac{\widehat{\nabla}\omega^0}{\Omega\omega^0}},
\frac{\widehat{\nabla}\omega^0}{\Omega\omega^0}\right)>0.
\end{equation}}
  \subsection{On trapped surfaces}
  Since the definition of black hole given here is not a priori the standard one, it appears obvious to
  question if there is a formation of trapped surfaces in this
  setting. We therefore presents here what can be viewed as
  an attempt to answer the question since the relation between
   the hypotheses formulated below and Cauchy data
  remains another important task to investigate.\\
  We first provide conditions which guarantee the existence
  of some particular codimension $2$ spacelike hypersurfaces in the region $\{m\geq
  0\}$. Namely, we assume that the vector field $\widehat{\nabla} m $ is timelike in                                                       $\{m\geq
  0\}$, this induces, under the condition $\frac{\partial m}{\partial\omega^0}<0$ that
  \begin{equation}\label{27}
\widehat{h}^{0i}\frac{\partial m}{\partial
  \omega^i}>0,\;\widehat{h}^{ij}\frac{\partial m}{\partial \omega^i}\frac{\partial
  m}{\partial\omega^j}>0;
  \end{equation}and
  \begin{equation}\label{28}
\left(\widehat{h}^{0i}\frac{\partial m}{\partial
  \omega^i}\right)^2-\widehat{h}^{00}\widehat{h}^{ij}\frac{\partial m}{\partial \omega^i}\frac{\partial
  m}{\partial\omega^j}>0;
  \end{equation}the hypersurfaces of constant $m$, $(\{m=\nu,\;\nu>0\})$ are spacelike. On the other hand, in $m>0$ the hypersurfaces of constant $\omega^0+m$,
   $(\{\omega^0+m=\mu\})$, are also spacelike.
   We consider the codimension $2$ spacelike surfaces in  $\{m>0\}$:
   \begin{equation*}
        H_{\mu\nu}:=\{(\omega^\alpha)\in
        \widehat{V},\;\omega^0+m=\mu,\;m=\nu,\;\mu>0,\;\nu>0\};
   \end{equation*}and suppose that they have a compact feeling (this is possible if $\{m=\nu,\;\nu>0\}$ is
   a Cauchy hypersurface and $H_{\mu\nu}$ is the boundary of a compact region of it ). At every point $\omega\in H_{\mu\nu}$, $(T_\omega
   H_{\mu\nu})^\perp$ has dimension $2$ and is timelike. Such a two plane
   cuts the future light cone at $\omega$ along two future null
   directions. We choose two future null vectors $l^+,\;l^-$
   supported by these directions up to a scaling by a positive
   factor --- $h(l^+,l^-)<0$. One can choose the scaling
   factor $\textbf{a}>0$ such that $(h(l^+,l^-) =-2)$. Taking $(e_a)$ a reference frame on $H_{\mu\nu}$, hence
   orthogonal to $l^+$ and $l^-$, we define the two
   null fundamental forms of $H_{\mu\nu}$ by:
   \begin{equation}\label{29}
\chi_{ab}=h(\widehat{\nabla}_{e_a}^{l^+},e_b),\;\underline{\chi}_{ab}=
h(\widehat{\nabla}_{e_a}^{l^-},e_b).
   \end{equation}The null mean curvatures of $H_{\mu\nu}$
   denoted $tr\chi$ and $tr\underline{\chi}$ are defined by:
\begin{equation}\label{30}
tr\chi=\widehat{h}^{ab}\chi_{ab},\;tr\underline{\chi}=\widehat{h}^{ab}\underline{\chi}_{ab}.
\end{equation}Now, at every point $\omega \in H_{\mu\nu}$,
the vector fields $\widehat{\nabla}\omega^0$ and $\widehat{\nabla}m$
belong to $(T_\omega
   H_{\mu\nu})^\perp$, and we assume $\widehat{\nabla}m$ to be future pointing, consequently
 \begin{equation}\label{31}
l^+=\textbf{a}(\textbf{b}^+\widehat{\nabla}\omega^0+\widehat{\nabla}m),\;
l^-=\textbf{a}(\textbf{b}^-\widehat{\nabla}\omega^0+\widehat{\nabla}m);
\end{equation}
\begin{equation*}
    \textbf{b}^{\pm}=\frac{-h(\widehat{\nabla}\omega^0,\widehat{\nabla}m)\pm
    \sqrt{(h(\widehat{\nabla}\omega^0,\widehat{\nabla}m))^2-
    h(\widehat{\nabla}\omega^0,\widehat{\nabla}\omega^0)h(\widehat{\nabla}m,\widehat{\nabla}m)}}
    {h(\widehat{\nabla}\omega^0,\widehat{\nabla}\omega^0)}.
\end{equation*}The two null forms read in terms of the components
 in the basis $(e_a)$ of
 the covariant derivatives (w.r.t. $e_a$) of the vector
fields $\widehat{\nabla}\omega^0$ and $\widehat{\nabla}m$ as:
\begin{equation}\label{32}
    \chi_{ab}=\textbf{a}\textbf{b}^+h(\widehat{\nabla}_{e_a}^{\widehat{\nabla}\omega^0},e_b)+
h(\widehat{\nabla}_{e_a}^{\widehat{\nabla}m},e_b)\equiv
\textbf{ab}^+\textbf{U}_{ab}+\textbf{a}\textbf{V}_{ab},
\end{equation}
\begin{equation}\label{33}
\underline{\chi}_{ab}=\textbf{a}\textbf{b}^-h(\widehat{\nabla}_{e_a}^{\widehat{\nabla}\omega^0},e_b)+
h(\widehat{\nabla}_{e_a}^{\widehat{\nabla}m},e_b)\equiv
\textbf{ab}^-\textbf{U}_{ab}+\textbf{a}\textbf{V}_{ab}.
\end{equation}
\subsection{On hypotheses on which trapped surfaces's formation is based}
From the expressions (\ref{33}) above it follows:
\begin{equation}\label{34}
    tr\chi=\textbf{ab}^+tr\textbf{U}+\textbf{a}tr\textbf{V},\;tr\underline{\chi}=\textbf{ab}^-tr\textbf{U}+\textbf{a}tr\textbf{V};
\end{equation}one deduces:
\begin{equation}\label{35}
    tr\textbf{V}=\frac{\textbf{b}^+tr\underline{\chi}-\textbf{b}^-tr\chi}{\textbf{a}(\textbf{b}^+-\textbf{b}^-)},
\end{equation}and hence to have
$tr\chi<0,\;tr\underline{\chi}<0$, it is necessary that $
tr\textbf{V}<0$. Further, the condition $tr\chi
tr\underline{\chi}>0$ requires that:
\begin{equation}\label{36}
    tr
    \textbf{U}\in\left]-\frac{tr\textbf{V}}{\textbf{b}^-},-\frac{tr\textbf{V}}{\textbf{b}^+}\right[.
\end{equation}Conversely if
\begin{equation}
  tr\textbf{V}<0,\;  tr
    \textbf{U}\in\left]-\frac{tr\textbf{V}}{\textbf{b}^-},-\frac{tr\textbf{V}}{\textbf{b}^+}\right[,
\end{equation}then
\begin{equation*}
    tr\chi <0,\;tr\underline{\chi}<0.
\end{equation*}
This previous analysis and the results of corollary \ref{co1} and
proposition \ref{p3}, induce the following theorem and its
corollary.
\subsection{Theorem -- Corollary}

 {\theorem {Let $(V,g)$ be a global
in time
 spacetime solution of the Einstein field
 equations where $g$ is of the following form in coordinates $(x^\alpha)$:
\begin{equation}
    g=-\frac{|\overline{g}|}{\Omega^2}(dx^0)^2+\overline{g}_{ij}dx^idx^j,\;\overline{g}=(\overline{g}_{ij}),
\end{equation} $\Omega\equiv \Omega (x^i)$ is an arbitrary positive
given scalar density satisfying $\Omega\rightarrow 0$ as
$r=|(x^i)|\rightarrow +\infty$. \\ Let denote $h=\Omega^2 g$
the conformal metric, $G$ its inverse, and $\widehat{\nabla}$ its associated riemannian
connection. Let be specified the gauge $(\omega^\alpha)$ for the spacetime $(V,g)$, with $
\omega^0=e^{-\Omega x^0},\;\omega^1=\Omega,\;\omega^a=\arctan
(\Omega x^a)$.\\ One suppose that $\mathcal{H}=\{m\equiv G\left(\frac{d\omega^0}{\Omega\omega^0},\frac{d\omega^0}{\Omega\omega^0}\right)=0\}\neq \emptyset$ and defines a (smooth) hypersurface, and the attached
 geometric objects (in $(\widehat{V},h)$), $\textbf{U},\textbf{V}$ and quantities $\textbf{a},\textbf{b}^\pm$ as described in
       (\ref{31})-(\ref{33}) satisfy:
  \begin{equation}
    h\left(\widehat{\nabla}_{\frac{\partial}{\partial\omega^0}}^{\frac{\widehat{\nabla}\omega^0}{\Omega\omega^0}},
\frac{\widehat{\nabla}\omega^0}{\Omega\omega^0}\right) <0, \;
   h\left(\widehat{\nabla}_{\widehat{h}^{0i}\frac{\partial}{\partial\omega^i}}^{\frac{\widehat{\nabla}\omega^0}{\Omega\omega^0}},
\frac{\widehat{\nabla}\omega^0}{\Omega\omega^0}\right)>0;
 \end{equation}
 \begin{equation}
  tr\textbf{V}<0,\;
      tr\textbf{U}\in\left]-\frac{tr\textbf{V}}{\textbf{b}^-}
      ,-\frac{tr\textbf{V}}{\textbf{b}^+}\right[.
       \end{equation}Then the spacetime $(V,g)$ undergoes trapped surfaces's formation.  }\label{th}}
{\corollary {One supposes that a gauge $(\omega^\alpha)$ is
specified as above for a spacetime $(V,g)$, and the extrinsic
curvature $K$ of the spacetime together with the attached
 geometric objects (in $(\widehat{V},h)$), $\textbf{U},\textbf{V}$ and quantities $\textbf{a},\textbf{b}^\pm$ as described in
       (\ref{31})-(\ref{33}) satisfy:
 \begin{equation}
   trK<0,\;K^{ij}\frac{\partial\omega^1}{\partial x^i}\frac{\partial\omega^1}{\partial
   x^j}\geq 0 ,\;
   h\left(\widehat{\nabla}_{\widehat{h}^{0i}\frac{\partial}{\partial\omega^i}}^{\frac{\widehat{\nabla}\omega^0}{\omega^1\omega^0}},
\frac{\widehat{\nabla}\omega^0}{\omega^1\omega^0}\right)>0;
 \end{equation}
 \begin{equation}
  tr\textbf{V}<0,\;
      tr\textbf{U}\in\left]-\frac{tr\textbf{V}}{\textbf{b}^-},
      -\frac{tr\textbf{V}}{\textbf{b}^+}\right[.
       \end{equation}Then the spacetime $(V,g)$ undergoes trapped surfaces's formation.  }}
       {\remark{To appreciate the physical interpretation of
       the various hypotheses above, one should first remember that $\widehat{\nabla}\omega^0$
        is viewed as time effect or time-flux, $m$ is the relative length of such flux,
         $\widehat{\nabla}m$} measures variation of such length.
          Furtheremore, the conditions required in the statements of the results above rest
           on components of the covariant derivatives of $\widehat{\nabla}\omega^0$ and $\widehat{\nabla}m$ in specific directions. }
\section{A lower bound problem for the total mass
 as a minimization problem via a Lagrangian formulation}
 The concept of mass or quasi-local mass is very important in general relativity.
 In [S.-T. Yau, Seminar on Differential Geometry (1982)], Penrose
 listed the search for a definition of such
 quasi-local mass as his number one problem in classical general
 relativity. Clearly, many important statements in general
 relativity require a good definition of quasi-local mass. This
 latter might help to control the dynamics of the gravitational
 field. It might also be used for energy methods for the analysis of
  hyperbolic equations
 in spacetimes. The positivity of the total mass is a matter of particular interest
as soon as such a quantity is defined
\cite{1},\cite{5}-\cite{9},\cite{40},\cite{47},\cite{48},
\cite{55},\cite{57},\cite{60},\cite{70}-\cite{72},\cite{74},\cite{79},\cite{90}.
This latter would be obtained if a positive lower bound inequality for the
total mass is
established. Let's recall that a lower bound for the total mass of
isolated systems could also guarantee the stability of such systems.
 In this section, according to the mass function adopted in this
 work, we analyze such a problem as a
minimization one under some prospective isoperimetric constraints
  whose nature (energy, area or volume,...) should depend on the desirable ends.
   For this purpose we use a foliation of the
conformal spacetime by the deformations of the black hole's event
horizon as follow. First of all let's remark that the actual event
horizon $\mathcal{H}^+$ is defined by an equation of the form
$\omega^0=X(\omega^1,...,\omega^n)$, and has an obvious
parametrization
\begin{equation*}
  H:\;  (\omega^1,\omega^a)\rightarrow
    (X(\omega^1,...,\omega^n),\omega^1,...,\omega^n).
\end{equation*}The Black hole's event horizon correspondingly is defined
by an equation of the form
\begin{equation*}
   \partial \mathcal{\textbf{B}}:\; \omega^0=\widetilde{X}(\omega^i);
\end{equation*}where
\begin{equation}
  \;(X(\omega^1,...,\omega^n),\omega^1,...,\omega^n)=
   (\widetilde{X}(\omega^1,...,\omega^n),\omega^1,...,\omega^n) \;\hbox{on $\mathcal{H}^+$}.
\end{equation}
 We denote by $\mathcal{H}_s,\;s\in
[-\epsilon,+\epsilon]$,
 a deformation of the "Black hole's event horizon"
$ \partial \mathcal{\textbf{B}}$ such that $\mathcal{H}_0\equiv
\partial \mathcal{\textbf{B}}$, adopted as follows:
\begin{equation*}
\mathcal{H}_0 \ni (\omega^\alpha)\mapsto
(\widetilde{T}(s)(\omega^i),\omega^j)\equiv ((T\circ
s(\omega^\alpha),\omega^1,...,\omega^n))\in \mathcal{H}_s.
\end{equation*}
Let denote $\underline{H}={]0,\Omega_{\max}]}\times
(]-\frac{\pi}{2},\;+\frac{\pi}{2}[)^{n-1}$, one has for $s\in
[-\epsilon,+\epsilon]$:
\begin{equation}
      \begin{array}{ll}
\mathcal{H}_s=\{(\omega^0,\omega^i)\in
\widehat{V}/\omega^0=\widetilde{T}(s)(\omega^i),\;(\omega^i)\in\underline{H}\},& \hbox{} \\
        \mathcal{H}_0=\{(\omega^0,\omega^i)\in \widehat{V}/\omega^0=
        \widetilde{X}(\omega^i),\;(\omega^i)\in\underline{H}\}. & \hbox{}
      \end{array}
\end{equation}The total mass ( squared) $M=\mathfrak{M}^2 $ of the spacetime is written
according to the deformation of the Black hole's event horizon as:
\begin{equation}
    M=\int_{-\epsilon}^{+\epsilon}\int_{\underline{H}}\frac{m^\ast(\widetilde{T}( s)(\omega^i),\omega^i)}
    {\widetilde{T}(s)}\frac{d \widetilde{T}(s)}{d s}
    \frac{d\omega^1}{\omega^1}d\omega^2...d\omega^n\;ds,
\end{equation}with $m^\ast=\omega^0\omega^1 m$. We remark that
$m^\ast$ is integrable up to the boundary. \\
Now set $\widetilde{T}(s)=T\circ s (.)$, then $\widetilde{T}(s)$ is
a distribution on $\underline{H}$ --- and define the functionals
$J$:
\begin{equation}
      J(\widetilde{T})=\int_{-\epsilon}^{+\epsilon}\int_{\underline{H}}
    \frac{m^\ast(\widetilde{T}(s)(\omega^i),\omega^j)}
    {\widetilde{T}(s)(\omega^i)}\frac{d \widetilde{T}(s)(\omega^i)}{ds}
    \frac{d\omega^1}{\omega^1}d\omega^2...d\omega^n\;ds,
\end{equation}which can then be rewritten as
\begin{equation}\label{e3}
  J(\widetilde{T}) =\int_{-\epsilon}^{+\epsilon}
    P(s,\widetilde{T}(s), \frac{d \widetilde{T}(s)}{ds})ds;
\end{equation}and and associate to $J$ another functional $I$:
\begin{equation}\label{e4}
    I(\widetilde{T})=\int_{-\epsilon}^{+\epsilon}
    P_0(s,\widetilde{T}(s), \frac{d \widetilde{T}(s)}{ds})ds;
    \end{equation}where the exact form of $P_0$ depends on the specific needed end.
This induces a minimization problem:
\begin{equation}
(P)\left\{
          \begin{array}{ll}
      \hbox{Find}\; \displaystyle\inf_{\widetilde{T}\in H^1(]-\epsilon;+\epsilon[)} J(\widetilde{T}), & \hbox{} \\
        I(\widetilde{T})=K, & \hbox{K is a constant.}
      \end{array}
      \right.
    \end{equation}where $H^1(]-\epsilon;+\epsilon[)$ is the classical Sobolev
    space. As one can notice we have given a Lagrangian formulation
    to the lower bound problem for the here considered total mass where the Lagrangian is the function
    $P(.,\widetilde{T}(.), \frac{d \widetilde{T}(.)}{ds})$.
    {\remark{In the geometric flow technique here, we have not considered the flow of the metric since at this stage we are not interested in a specific end for the spacetime. The natural metric on any subspace here is the induced one.}}
    \subsection{Frechet derivative $J'(\widetilde{T})$}
    The function $P(.,.,.)$ is smooth and we find the
    Frechet-derivative of $J$. The
    Isoperimetric constraint $I(\widetilde{T})=K$ is combined with
    the boundary conditions that induce the set of constraints
    (we search for distribution $\widetilde{T}\in
    H^1(]-\epsilon;+\epsilon[)$)
    \begin{equation}
      C=\{  \widetilde{T}\in
    H^1(]-\epsilon;+\epsilon[),\; \widetilde{T}(+\epsilon)\equiv 0,\;\widetilde{T}(-\epsilon)\equiv
        \widetilde{Z} \}.
    \end{equation} $\widetilde{Z}$ is such that the $0$-set of
     $\omega^0-\widetilde{Z}$ defines the initial hypersurface. Obviously $\widetilde{Z}$ is given by:
    \begin{equation*}
                  \begin{array}{ll}
          \widetilde{Z}: \widehat{V}\rightarrow \mathbb{R} & \hbox{} \\
            (\omega^0,\omega^i)\mapsto \widetilde{Z}(\omega^0,\omega^i)=e^{-\omega^1Z(\omega^i)},\;
            Z(\omega^i)\equiv x^0(\omega^1,...,\omega^n). & \hbox{}
          \end{array}
           \end{equation*}For $u\in C$, the set of
    admissible directions of $u$ in the sense of Frechet denoted $K(u)$ is:
    \begin{equation}
          \begin{array}{ll}
        K(u)=:\{w\in
    H^1(]-\epsilon;+\epsilon[)/  \;\hbox{there exists a sequence}\; (w_n)\in
       H^1(]-\epsilon;+\epsilon[)\; \hbox{}  & \hbox{} \\
      \hbox{converging to $w$  and a sequence}\; (e_n),\;e_n > 0 \; \hbox{such that}\; u+e_n w_n\in C \}.& \hbox{}
      \end{array}
    \end{equation}If $u+e_n w_n\in C$ then one can deduce that
    $w_n(+\epsilon)=0,\;w_n(-\epsilon)=0$, and since functions of $H^1(]-\epsilon;+\epsilon[)
    $ are continuous at the boundary, and the application trace is
    also continuous, one deduces that
    $w(-\epsilon)=0,\;w(+\epsilon)=0$. Conversely if
    $w(-\epsilon)=0,\;w(+\epsilon)=0$, one constructs
    $u+\frac{1}{n}w$ which satisfies the constraints. It follows that
    \begin{equation}
K(u)=H_0^1(]-\epsilon;+\epsilon[).
    \end{equation}
    For $w\in H_0^1(]-\epsilon;+\epsilon[) $,
    \begin{equation}\label{}
        (J'(u),w)=\lim_{\lambda\rightarrow 0}\frac{J(u+\lambda w)-J(u)}{\lambda}
    \end{equation}
\begin{equation}
= \lim_{\lambda\rightarrow
0}\frac{1}{\lambda}\int_{-\epsilon}^{+\epsilon}\{P(u+\lambda
w,u'+\lambda w')-P(u,u')\}ds.
\end{equation}
Since the application $(\widetilde{T},\widetilde{T}')\mapsto
P((\widetilde{T},\widetilde{T}'))$ is differentiable, one can write:
\begin{equation*}
    P(u+\lambda
w,u'+\lambda w')=P(u,u')+DP(u,u'). \lambda (w,w')
+|\lambda|\|(w,w')\|O(\lambda (w,w')),
\end{equation*}where $DP$ denotes the differential of $P$ --- this implies that
\begin{equation*}
    (J'(u),w)=\int_{-\epsilon}^{+\epsilon}DP(u,u').(w,w') ds
    \end{equation*}
\begin{equation*}
    = \int_{-\epsilon}^{+\epsilon} \left[\frac{\partial P}{\partial \widetilde{T}}(u,u')w +
    \frac{\partial P}{\partial \widetilde{T}'}(u,u')w'\right] ds
     \end{equation*}
\begin{equation*}
    =\int_{-\epsilon}^{+\epsilon} \left[\frac{\partial P}{\partial \widetilde{T}}(u,u').w\right]ds+
   \int_{-\epsilon}^{+\epsilon} \left[\frac{d}{ds}\left( \frac{\partial P}{\partial \widetilde{T}'}(u,u').w\right)-
   \frac{d}{ds}\left(\frac{\partial P}{\partial \widetilde{T}'}(u,u')\right).w\right] ds
    \end{equation*}
\begin{equation*}
    =\int_{-\epsilon}^{+\epsilon} \left[\frac{\partial P}{\partial \widetilde{T}}(u,u')-
    \frac{d}{ds}\left( \frac{\partial P}{\partial \widetilde{T}'}(u,u')\right)\right]. w\;ds+
    \left[ \frac{\partial P}{\partial \widetilde{T}'}(u,u').w \right] _{-\epsilon}^{+\epsilon}
     \end{equation*}
\begin{equation*}
    = \int_{-\epsilon}^{+\epsilon} \left[\frac{\partial P}{\partial \widetilde{T}}(u,u')-
    \frac{d}{ds}\left( \frac{\partial P}{\partial \widetilde{T}'}(u,u')\right)\right]. w\;ds+
     \end{equation*}
\begin{equation*}
   \frac{\partial P}{\partial \widetilde{T}'}(u(+\epsilon),u'(+\epsilon))\delta_{+\epsilon}(w)-
    \frac{\partial P}{\partial
    \widetilde{T}'}(u(-\epsilon),u'(-\epsilon))\delta_{-\epsilon}(w);
     \end{equation*}
\begin{equation*}
\end{equation*}$\delta_{\pm\epsilon}$ denotes Dirac-distribution at $\pm\epsilon$,
one deduces that:
\begin{equation*}
   J'(u)= \frac{\partial P}{\partial \widetilde{T}}(u,u')-
    \frac{d}{ds}\left( \frac{\partial P}{\partial \widetilde{T}'}(u,u')\right)+
     \end{equation*}
\begin{equation}
   \frac{\partial P}{\partial \widetilde{T}'}(u(+\epsilon),u'(+\epsilon))\delta_{+\epsilon}-
    \frac{\partial P}{\partial
    \widetilde{T}'}(u(-\epsilon),u'(-\epsilon))\delta_{-\epsilon}.
     \end{equation}
     Let's recall that:
\begin{equation*}
   P(\widetilde{T},\frac{d\widetilde{T}}{ds})= \int_{\underline{H}}
    \frac{m^\ast(\widetilde{T}(s)(\omega^i),\omega^i)}
    {\widetilde{T}(s)(\omega^i)}\frac{d \widetilde{T}(s)(\omega^i)}{ds}
    \frac{d\omega^1}{\omega^1}d\omega^2...d\omega^n;
        \end{equation*}
therefore
    \begin{equation}
        \frac{\partial P}{\partial \widetilde{T}}=\int_{\underline{H}}
        \frac{\widetilde{T}(s)(\omega^i)\frac{\partial m^\ast}{\partial \omega^0}(\widetilde{T}(s)(\omega^i))-
        m^\ast(\widetilde{T}(s)(\omega^i),\omega^i)}{\widetilde{T}^2(s)(\omega^i)}\frac{d
        \widetilde{T}(s)}{ds}(\omega^i)\frac{d\omega^1}{\omega^1}d\omega^2...d\omega^n,
    \end{equation}
    \begin{equation}
        \frac{\partial P}{\partial \widetilde{T}'}=\int_{\underline{H}}
        \frac{m^\ast(\widetilde{T}(s)(\omega^i),\omega^j)}{\widetilde{T}(s)(\omega^i)}
        \frac{d\omega^1}{\omega^1}d\omega^2...d\omega^n,
    \end{equation}
\begin{equation*}
        \frac{d}{ds}\left(\frac{\partial P}{\partial
    \widetilde{T}'}\right)=\int_{\underline{H}}\frac{d}{ds}
    \left(\frac{m^\ast(\widetilde{T}(s)(\omega^i),\omega^i)}{\widetilde{T}(s)(\omega^i)}\right)\frac{d\omega^1}
    {\omega^1}d\omega^2...d\omega^n
\end{equation*}
\begin{equation}
=\int_{\underline{H}}
        \frac{\widetilde{T}(s)(\omega^i)\frac{\partial m^\ast}{\partial \omega^0}(\widetilde{T}(s)(\omega^i),\omega^i)-
        m^\ast(\widetilde{T}(s)(\omega^i),\omega^i)}{\widetilde{T}^2(s)(\omega^i)}\frac{d
        \widetilde{T}(s)}{ds}\frac{d\omega^1}{\omega^1}d\omega^2...d\omega^n.
\end{equation}
It appears that $P$ satisfies the Lagrange's equation
\begin{equation}
\frac{\partial P}{\partial
\widetilde{T}}-\frac{d}{ds}\left(\frac{\partial P}{\partial
    \widetilde{T}'}\right)=0,\;\forall \widetilde{T}.
\end{equation}
\subsection{Euler's condition for $R\equiv P+\lambda P_0$ -- A fundamental equation}
 We recall here that if $J(\widetilde{T})$ is an extremum for the functional $J(.)$ subject to the isoperimetric constraint $I(.)=K$, then $\widetilde{T}$ is an extremal of $\int R (.)=\int (P+\lambda
P_0) (.)$ for some constant $\lambda $. The initially unknown multiplier $\lambda$ must be determined at the end using the isoperimetric constraint
$\int P_0=K$\\  The Euler's equation for $R=P+\lambda
P_0$ reads:
\begin{equation*}
    \frac{\partial P}{\partial
\widetilde{T}}-\frac{d}{ds}\left(\frac{\partial P}{\partial
    \widetilde{T}'}\right)+\lambda \left[ \frac{\partial P_0}{\partial
\widetilde{T}}-\frac{d}{ds}\left(\frac{\partial P_0}{\partial
    \widetilde{T}'}\right)\right]+
    \end{equation*}
\begin{equation*}
    \frac{\partial P}{\partial \widetilde{T}'}(\widetilde{T}(+\epsilon),\widetilde{T}'(+\epsilon))\delta_{+\epsilon}-
    \frac{\partial P}{\partial
    \widetilde{T}'}(\widetilde{T}(-\epsilon),\widetilde{T}'(-\epsilon))\delta_{-\epsilon}+
     \end{equation*}
\begin{equation}
    \lambda\frac{\partial P_0}{\partial \widetilde{T}'}(\widetilde{T}(+\epsilon),\widetilde{T}'(+\epsilon))\delta_{+\epsilon}-
    \lambda\frac{\partial P_0}{\partial
    \widetilde{T}'}(\widetilde{T}(-\epsilon),\widetilde{T}'(-\epsilon))\delta_{-\epsilon}=0.
\end{equation}
Considering that $P$ satisfies the equation $\frac{\partial P}{\partial
\widetilde{T}}-\frac{d}{ds}\left(\frac{\partial P}{\partial
    \widetilde{T}'}\right)=0,\;\forall \widetilde{T}$, this latter equation
reduces then for $\lambda\neq 0$ to the following that we consider as a fundamental equation:
\begin{equation*}
    \frac{\partial P_0}{\partial
\widetilde{T}}-\frac{d}{ds}\left(\frac{\partial P_0}{\partial
    \widetilde{T}'}\right)+\left[\frac{1}{\lambda}\frac{\partial P}{\partial \widetilde{T}'}(\widetilde{T}(+\epsilon),\widetilde{T}'(+\epsilon))
+ \frac{\partial P_0}{\partial \widetilde{T}'}(\widetilde{T}(+\epsilon),\widetilde{T}'(+\epsilon))\right]\delta_{+\epsilon}-
    \end{equation*}
\begin{equation}\label{e1}
\left[\frac{1}{\lambda}\frac{\partial P}{\partial \widetilde{T}'}(\widetilde{T}(-\epsilon),\widetilde{T}'(-\epsilon))
+ \frac{\partial P_0}{\partial \widetilde{T}'}(\widetilde{T}(-\epsilon),\widetilde{T}'(-\epsilon))\right]\delta_{-\epsilon}
   =0.
\end{equation}
\subsection{Example of $P_0$-- A Sturn- Liouville operator}
Here, we consider a situation where the perturbation of the black hole is due to a
 forced term $f$ and the interested functional is given by
  \begin{equation*}
 I(\widetilde{T})=   \int_{-\epsilon}^{+\epsilon}P_0(s,\widetilde{T}(s),\widetilde{T}'(s))ds=
 \int_{-\epsilon}^{+\epsilon}\left(\frac{1}{2}(\widetilde{T}'(s))^2+\frac{1}{2}(\widetilde{T}(s))^2-f(s)\widetilde{T}(s)\right)ds.
  \end{equation*}A minimizer of $J(\widetilde{T})$ subject to the isoperimetric constraint $ I(\widetilde{T})= K$ must be a function $\widetilde{T}_\ast (s,\lambda)$ satisfying the differential equation
    \begin{equation*}
    -\widetilde{T}^{''}+\widetilde{T}-f +\left[\frac{1}{\lambda}\frac{\partial P}{\partial \widetilde{T}'}(\widetilde{T}(+\epsilon),\widetilde{T}'(+\epsilon))
+ \frac{\partial P_0}{\partial \widetilde{T}'}(\widetilde{T}(+\epsilon),\widetilde{T}'(+\epsilon))\right]\delta_{+\epsilon}-
    \end{equation*}
\begin{equation}\label{e5}
\left[\frac{1}{\lambda}\frac{\partial P}{\partial \widetilde{T}'}(\widetilde{T}(-\epsilon),\widetilde{T}'(-\epsilon))
+ \frac{\partial P_0}{\partial \widetilde{T}'}(\widetilde{T}(-\epsilon),\widetilde{T}'(-\epsilon))\right]\delta_{-\epsilon}
   =0,
\end{equation}and where the parameter $\lambda$ enjoys $ I(\widetilde{T}_\ast)= K$. This equation presents a Sturn-Liouville operator and one might expect quantization in its resolution process.
\subsection{Second derivative $J"(\widetilde{T})\widetilde{V}$}
In order to characterize the minimality condition we are now
studying
\begin{equation*}
    \frac{J'(\widetilde{T}+\sigma \widetilde{V},w)-J'(\widetilde{T},w)}{\sigma}.
\end{equation*}Let's recall that:
\begin{equation}\label{}
\begin{array}{ll}
   J'(\widetilde{T})= \underset{=0}{\underbrace{\frac{\partial P}{\partial \widetilde{T}}
   (\widetilde{T},\widetilde{T}')-
    \frac{d}{ds}\left( \frac{\partial P}{\partial \widetilde{T}'}(\widetilde{T},\widetilde{T}')\right)}}+& \hbox{} \\
   \frac{\partial P}{\partial \widetilde{T}'}(\widetilde{T}(+\epsilon),\widetilde{T}'(+\epsilon))\delta_{+\epsilon}-
    \frac{\partial P}{\partial
    \widetilde{T}'}(\widetilde{T}(-\epsilon),\widetilde{T}'(-\epsilon))\delta_{-\epsilon}.
     & \hbox{}
  \end{array}
\end{equation}We denote by $S$ the application
\begin{equation}\label{9}
    S: (\widetilde{T},\widetilde{T}')\rightarrow \frac{\partial P}{\partial
\widetilde{T}'}=\int_{\underline{H}}\frac{m^\ast(\widetilde{T})}{\widetilde{T}}\frac{d\omega^1}{\omega^1}
d\omega^2 ...\omega^n,
\end{equation}which is independent of $\widetilde{T}'$ and is
differentiable for every $(\widetilde{T},\widetilde{T}')$.
\begin{equation*}
    J'(\widetilde{T}+\sigma \widetilde{V},w)=
\end{equation*}
\begin{equation*}
  \left\langle\int_{\underline{H}}\frac{m^\ast(\widetilde{T}
    (+\epsilon)+\sigma \widetilde{V}(+\epsilon))}{\widetilde{T}(+\epsilon)+
    \sigma\widetilde{V}(+\epsilon)}\frac{d\omega^1}{\omega^1}d\omega^2
    ...d\omega^n\delta_{+\epsilon},w\right\rangle
    -
\end{equation*}
\begin{equation*}
    \left\langle\int_{\underline{H}}\frac{m^\ast(\widetilde{T}
    (-\epsilon)+\sigma \widetilde{V}(-\epsilon))}{\widetilde{T}(-\epsilon)+
    \sigma\widetilde{V}(-\epsilon)}\frac{d\omega^1}{\omega^1}d\omega^2
    ...d\omega^n\delta_{-\epsilon},w \right\rangle
    .
\end{equation*}One has $\frac{\partial S}{\partial \widetilde{T}'}\equiv
0$ and:
\begin{equation*}
    \frac{\partial S}{\partial \widetilde{T}}=\int_{\underline{H}}\frac{\frac{\partial
    m^\ast}{\partial\omega^0}(\widetilde{T})\times
    \widetilde{T}-m^\ast(\widetilde{T})}{\widetilde{T}^2}\frac{d\omega^1}{\omega^1}...d\omega^n,
\end{equation*}and
\begin{equation*}
    DS(\widetilde{T},\widetilde{T}')(U,W)=\frac{\partial S}{\partial
    \widetilde{T}}. U+0. W\equiv \frac{\partial S}{\partial
    \widetilde{T}}.U.
\end{equation*}One can then write:
\begin{equation*}
    \left\langle S\left((\widetilde{T}(+\epsilon)+\sigma \widetilde{V}(+\epsilon),
    \widetilde{T}'(+\epsilon)+\sigma \widetilde{V}'(+\epsilon))\right)\delta_{+\epsilon},w \right\rangle
\end{equation*}
\begin{equation*}
    =\left\langle \left\{S(\widetilde{T}(\epsilon),
    \widetilde{T}'(\epsilon))+\sigma DS (\widetilde{T}(\epsilon),
    \widetilde{T}'(\epsilon))(\widetilde{V}(\epsilon),\widetilde{V}'(\epsilon))\right\}\delta_{+\epsilon},w
    \right\rangle
\end{equation*}
\begin{equation*}
    + \left\langle \sigma O(\sigma),w \right\rangle,
\end{equation*}The same holds substituting $+\epsilon$ by $-\epsilon$.
One then deduces that:
\begin{equation*}
    \frac{\left\langle J'(\widetilde{T}+\sigma \widetilde{V}),w\right\rangle-
    \left\langle J'(\widetilde{T}),w\right\rangle}{\sigma}=
\end{equation*}
\begin{equation*}
    \left\langle \left\{DS (\widetilde{T}(\epsilon),
    \widetilde{T}'(\epsilon))(\widetilde{V}(\epsilon),\widetilde{V}'(\epsilon))\right\}
    \delta_{+\epsilon},w\right\rangle-
\end{equation*}
\begin{equation*}
    \left\langle \left\{DS (\widetilde{T}(-\epsilon),
    \widetilde{T}'(-\epsilon))(\widetilde{V}(-\epsilon),\widetilde{V}'(-\epsilon))\right\}\delta_{-\epsilon},w
    \right\rangle
       +\left\langle O(\sigma),w \right\rangle,
          \end{equation*}and
\begin{equation*}
  \left\langle   J"(\widetilde{T})\widetilde{V},w\right\rangle=
\left\langle \left\{\frac{\partial S}{\partial \widetilde{T}}
(\widetilde{T}(+\epsilon),
\widetilde{T}'(+\epsilon)).\widetilde{V}(+\epsilon)\right\}\delta_{+\epsilon},w\right\rangle
\end{equation*}
\begin{equation*}
-\left\langle \left\{\frac{\partial S}{\partial \widetilde{T}}
(\widetilde{T}(-\epsilon),
\widetilde{T}'(-\epsilon)).\widetilde{V}(-\epsilon)\right\}\delta_{-\epsilon},w\right\rangle;
\end{equation*}using the expression of $\frac{\partial S}{\partial
\widetilde{T}}$ above, one has finally:
\begin{equation*}
    \left\langle   J"(\widetilde{T})\widetilde{V},w\right\rangle  =
\end{equation*}
\begin{equation*}
 \left\langle
   \int_{\underline{H}} \frac{\frac{\partial m^\ast}{\partial\omega^0}(\widetilde{T}(+\epsilon)).
    \widetilde{T}(+\epsilon)-m^\ast(\widetilde{T}(+\epsilon))}{\widetilde{T}^2(+\epsilon)}
   \frac{d\omega^1}{\omega^1}d\varpi \;\widetilde{V}(+\epsilon)\delta_{+\epsilon},w\right\rangle
\end{equation*}
\begin{equation}
    -\left\langle \int_{\underline{H}}
    \frac{\frac{\partial m^\ast}{\partial\omega^0}(\widetilde{T}(-\epsilon)).
    \widetilde{T}(-\epsilon)-m^\ast(\widetilde{T}(-\epsilon))}{\widetilde{T}^2(-\epsilon)}
    \frac{d\omega^1}{\omega^1}d\varpi\;\widetilde{V}(-\epsilon)\delta_{-\epsilon},w\right\rangle.
\end{equation}
\subsection{On a second order condition for a minimum}
As a result of the expression of $\left\langle
J"(\widetilde{T})\widetilde{V},w\right\rangle $ above, it follows
that:
\begin{equation*}
    \left\langle   J"(\widetilde{T})\widetilde{V},\widetilde{V}\right\rangle  =
\end{equation*}
\begin{equation*}
   \int_{\underline{H}} \frac{\frac{\partial m^\ast}{\partial\omega^0}(\widetilde{T}(+\epsilon)).
    \widetilde{T}(+\epsilon)-m^\ast(\widetilde{T}(+\epsilon))}{\widetilde{T}^2(+\epsilon)}
   \frac{d\omega^1}{\omega^1}d\varpi\;(\widetilde{V}(+\epsilon))^2
\end{equation*}
\begin{equation*}
    -\int_{\underline{H}}
    \frac{\frac{\partial m^\ast}{\partial\omega^0}(\widetilde{T}(-\epsilon)).
    \widetilde{T}(-\epsilon)-m^\ast(\widetilde{T}(-\epsilon))}{\widetilde{T}^2(-\epsilon)}
    \frac{d\omega^1}{\omega^1}d\varpi\;(\widetilde{V}(-\epsilon))^2.
\end{equation*}Since $\widetilde{V}(+\epsilon)$ in any case is required to
satisfy $\widetilde{V}(+\epsilon)=0$, it follows that the consistent
term is
\begin{equation}
 \left\langle   J"(\widetilde{T})\widetilde{V},\widetilde{V}\right\rangle  =
 -\int_{\underline{H}}
    \frac{\frac{\partial m^\ast}{\partial\omega^0}(\widetilde{T}(-\epsilon)).
    \widetilde{T}(-\epsilon)-m^\ast(\widetilde{T}(-\epsilon))}{\widetilde{T}^2(-\epsilon)}
    \frac{d\omega^1}{\omega^1}d\varpi\;(\widetilde{V}(-\epsilon))^2.
\end{equation}Now one wants to know if
\begin{equation*}
    \left\langle
    J"(\widetilde{T}_\ast)\widetilde{V},\widetilde{V}\right\rangle\;\geq
    0,
\end{equation*}independently of the value of
$\widetilde{V}(-\epsilon)$. It appears therefore necessary to study
the term
 \begin{equation*}
    A(\widetilde{T}(s),\omega^i)=\frac{\partial m^\ast}{\partial\omega^0}(\widetilde{T}(s)(\omega^i),\omega^i).
    \widetilde{T}(s)(\omega^i)-m^\ast(\widetilde{T}(s)(\omega^i),\omega^i)
\end{equation*}since $\left\langle J"(\widetilde{T})\widetilde{V},\widetilde{V}\right\rangle \geq 0$ if
\begin{equation}
 A(\widetilde{T}(-\epsilon),\omega^i)\leq 0.
\end{equation}We proceed to an analysis of the function
$A(\omega^0,\omega^i)=\omega^0\frac{\partial
m^\ast}{\partial\omega^0}(\omega^0,\omega^i)-m^\ast(\omega^0,\omega^i)$.
First of all, we have that
\begin{equation*}
    \omega^0\frac{\partial
    m^\ast}{\partial\omega^0}-m^\ast=\omega^1(\omega^0)^2\frac{\partial
    m}{\partial\omega^0}.
\end{equation*}The sign of $A(\omega^0,\omega^i)$ is then the one of
$\frac{\partial
    m}{\partial\omega^0}$ whose expression is given by the relation
    (\ref{5}) above. It results the following proposition.
{\proposition{} The function
$-A(\omega^0,\omega^i)=m^\ast-\omega^0\frac{\partial
m^\ast}{\partial\omega^0}$ is non negative if the extrinsic
curvature $K$ of the spacetime satisfies the inequalities
\begin{equation}\label{6}
    trK< 0,\;K^{ij}\frac{\partial\omega^1}{\partial x^i}\frac{\partial\omega^1}{\partial x^j}\geq 0.
\end{equation}
 }
 \\
 \textbf{Proof}\\
 The proof is an obvious consequence of the sign of
 $-A(\omega^\alpha)\equiv -\omega^1\omega^0\frac{\partial m}{\partial\omega^0}$.\\
 This proposition induces that independently of the value of
 $\widetilde{V}(-\epsilon)$,
  $\left\langle   J"(\widetilde{T})\widetilde{V},\widetilde{V}\right\rangle \geq
  0$ under the hypotheses (\ref{6}).
\subsection{Theorem -- Corollary }
{\theorem {One supposes that a gauge $(\omega^\alpha)$ is specified as in theorem \ref{th} for a spacetime $(V,g)$, $\mathcal{H}=\{m\equiv G\left(\frac{d\omega^0}{\Omega\omega^0},\frac{d\omega^0}{\Omega\omega^0}\right)=0\}\neq \emptyset$ and defines a (smooth) hypersurface,
  and the attached
 geometric objects (in $(\widehat{V},h)$) $\textbf{U},\textbf{V}$, and quantities $\textbf{a},\textbf{b}^\pm$ as described in
       (\ref{31})-(\ref{33}) satisfy:
  \begin{equation}
    h\left(\widehat{\nabla}_{\frac{\partial}{\partial\omega^0}}^{\frac{\widehat{\nabla}\omega^0}{\Omega\omega^0}},
\frac{\widehat{\nabla}\omega^0}{\Omega\omega^0}\right) <0, \;
   h\left(\widehat{\nabla}_{\widehat{h}^{0i}\frac{\partial}{\partial\omega^i}}^{\frac{\widehat{\nabla}\omega^0}{\Omega\omega^0}},
\frac{\widehat{\nabla}\omega^0}{\Omega\omega^0}\right)>0;
 \end{equation}
 \begin{equation}
  tr\textbf{V}<0,\;
      tr\textbf{U}\in\left]-\frac{tr\textbf{V}}{\textbf{b}^-}
      ,-\frac{tr\textbf{V}}{\textbf{b}^+}\right[.
       \end{equation} \\One supposes further that for a given smooth function $P_0$  of three real variables, the fundamental equation (\ref{e1}) (with $P$ described in (\ref{e3}))
              \begin{equation*}
    \frac{\partial P_0}{\partial
\widetilde{T}}-\frac{d}{ds}\left(\frac{\partial P_0}{\partial
    \widetilde{T}'}\right)+\left[\frac{1}{\lambda}\frac{\partial P}{\partial \widetilde{T}'}(\widetilde{T}(+\epsilon),\widetilde{T}'(+\epsilon))
+ \frac{\partial P_0}{\partial \widetilde{T}'}(\widetilde{T}(+\epsilon),\widetilde{T}'(+\epsilon))\right]\delta_{+\epsilon}-
    \end{equation*}
\begin{equation}\label{e6}
\left[\frac{1}{\lambda}\frac{\partial P}{\partial \widetilde{T}'}(\widetilde{T}(-\epsilon),\widetilde{T}'(-\epsilon))
+ \frac{\partial P_0}{\partial \widetilde{T}'}(\widetilde{T}(-\epsilon),\widetilde{T}'(-\epsilon))\right]\delta_{-\epsilon}
   =0
\end{equation}
       admits a solution $\widetilde{T}_\ast (.,\lambda)$ with $\lambda$ satisfying $I(\widetilde{T}_\ast)=K$.\\ Then the total mass $\mathfrak{M}$ corresponding to
the mass function
$m=G(\frac{d\omega^0}{\omega^1\omega^0},\frac{d\omega^0}{\omega^1\omega^0})\equiv
h(\frac{\widehat{\nabla}\omega^0}{\omega^1\omega^0},\frac{\widehat{\nabla}\omega^0}{\omega^1\omega^0})$
(with $G$ the inverse of $h\equiv \Omega^2g $) admits a lower bound.  }\label{th1}}\\
\textbf{Proof}\\The proof is already done since $\widetilde{T}_\ast (.,\lambda)$ is an extremal of $J(\widetilde{T})$ and the assumptions in the theorem imply that $\left\langle   J"(\widetilde{T}_\ast)\widetilde{V},\widetilde{V}\right\rangle \geq
  0$ as established above.
{\corollary {One supposes that a gauge $(\omega^\alpha)$ is
specified as above for a spacetime $(V,g)$, and the extrinsic
curvature $K$ of the spacetime together with the attached
 geometric objects (in $(\widehat{V},h)$) $\textbf{U},\textbf{V}$, and quantities $\textbf{a},\textbf{b}^\pm$ as described in
       (\ref{31})-(\ref{33}) satisfy:
 \begin{equation}
   trK<0,\;K^{ij}\frac{\partial\omega^1}{\partial x^i}\frac{\partial\omega^1}{\partial
   x^j}\geq 0 ,\;
   h\left(\widehat{\nabla}_{\widehat{h}^{0i}\frac{\partial}{\partial\omega^i}}^{\frac{\widehat{\nabla}\omega^0}{\omega^1\omega^0}},
\frac{\widehat{\nabla}\omega^0}{\omega^1\omega^0}\right)>0;
 \end{equation}
 \begin{equation}
  tr\textbf{V}<0,\;
      tr\textbf{U}\in\left]-\frac{tr\textbf{V}}{\textbf{b}^-},
      -\frac{tr\textbf{V}}{\textbf{b}^+}\right[.
       \end{equation} \\ One supposes further that for a given smooth function $P_0$  of three real variables, the equation (\ref{e6})
       admits a solution $\widetilde{T}_\ast (.,\lambda)$ with $\lambda$ satisfying $I(\widetilde{T}_\ast)=K$. Then total mass $\mathfrak{M}$ corresponding to the mass function
$m=G(\frac{d\omega^0}{\omega^1\omega^0},\frac{d\omega^0}{\omega^1\omega^0})\equiv
h(\frac{\widehat{\nabla}\omega^0}{\omega^1\omega^0},\frac{\widehat{\nabla}\omega^0}{\omega^1\omega^0})$
 admits a lower bound. }}
\section{Conclusion and Outlook}
There are many theorems and conjectures about black holes difficult
to prove and even often to state in a precise way
\cite{11},\cite{17},\cite{24},\cite{31},\cite{36},\cite{37},\cite{42},
\cite{45},\cite{66},\cite{77}, \cite{86}. In this paper, we have
proposed a framework with the hope that some aspects of issues of singularities
and black holes might have a chance to be understood though the relation
between the current approach and standard ones should be examined
carefully. Indeed, considering the crucial and central question of
mass or quasi-local mass in general relativity, we have proposed a
new mass function and for this mass we have established a possibility of the existence of a
lower bound for the corresponding "total mass" for a system collapsing to a black
hole thanks to some
hypotheses on the extrinsic curvature of the spacetime. Such
inequality in general is likely to guarantee the stability of
isolated systems, furthermore, it can help as a tool in the analysis
of partial differential equations (energy estimates,...) in black
holes's geometry. Concerning the question of existence or the nature
of horizons, some information were gained through the analysis
above, however, a profound analysis of the relation between the
horizon as described here with the standard definition or the one
attributed to Hayward Sean A. \cite{97} would be interesting, probably in relation with the detection of horizons using curvature invariants. For
the question of singularities and related problems (the fate of the
Cauchy horizons when they exist,
  cosmic censorship conjectures,...
\cite{2},\cite{3},\cite{10}-\cite{16},\cite{22}-\cite{37},\cite{82}),
it should be analyzed according to the approach here in a subsequent
work starting from standard black holes. Importantly, the
characterization of trapped surfaces in this setting or analysis of
the topology of the obtained black hole in general should be done
\cite{17},\cite{77}. The framework of this paper offers also new
possibilities for the analysis of waves on fixed backgrounds
 \cite{2},\cite{3},\cite{10}-\cite{16},\cite{22}-\cite{39},\cite{41}-\cite{44},
\cite{46},\cite{48}-\cite{52},\cite{58}-\cite{61},\cite{82}-\cite{84}),
conformal scattering is one such possibilities. Let's recall that
the choice of an appropriate gauge is at the heart of measurement of
decays of waves and of the stability problem \cite{93},\cite{94}. On
the other hand, the conditions that guarantee the existence of black
holes's region in this context should be compared if possible to
other energy conditions involving the energy momentum tensor. \\All
the results here are based on a mass function $\textbf{m}$ whose
relations with other mass-functions
(\cite{20},\cite{21},\cite{79},\cite{85}-\cite{90},\cite{92}) if
established would provide new insights to the understanding of
implications of general relativity. This mass-function is not always
positive, but the positivity of the "total mass" of the system is
obtainable provided some conditions are given on the isoperimetric constraint;
 this offers really many perspectives as one can appreciate through the example (\ref{e5}).
  However, it is worth noting that there are unified theory
projects in the literature based on Einstein-Yang-Mills-Dirac equations where the
system has negative energy and hence does not satisfy the positivity
conditions in the Penrose-Hawking singularity theorem. \\The results
obtained here assume the existence of a global in time space-time,
it is essential to study the global solvability of the wave
equations and quasi-linear wave equations (Einstein equations)
together with the analysis of the properties of solutions and their
decays. As illustration, there is a well known great interest to
enquire what happens to solutions of waves equations when Cauchy
horizons occur. Novel approaches to the global study of nonlinear
hyperbolic equations based on microlocal analysis have been proved
successful in the case of cosmological black holes
\cite{43}-\cite{45}. At least the approach of this paper is based on
the conformal embedding of the physical spacetime into a manifold
with corners, this permits to envisage using of methods of
microlocal analysis of Hintz and Vasy \cite{43}-\cite{45} for
various problems under consideration.
%
\acknowledgements
The author thanks the organizers of the two trimester programs
(Mathematical General Relativity, MSRI, Berkeley-California 2013 --
Mathematical General Relativity, IHP, Paris 2015) for the financial
support to his participation to these programs, which has
comforted his interest in Mathematical general relativity.\\
 The author acknowledges also fruitful collaborations with colleagues of
the Research Unit of Mathematics and Applications (RUMA) in the
Department of Mathematics and Computer science of the University of
Dschang, West region, Cameroon.
%
 \section*{Conflict of interest}
 The author declares that he has no conflict of interest..
\bibliographystyle{}
\bibliography{Bibliography article JMP}

\end{document}